\title{Nonabelian Hodge theory in positive characteristic for rank-2 vector bundles over curves}

\documentclass[11pt]{amsart}
\usepackage{mathrsfs,eucal,amsmath,amssymb,latexsym,dsfont,stmaryrd}
\usepackage{enumerate}
\usepackage[square,numbers]{natbib}
\usepackage{charter}
\usepackage{eulervm} 
\usepackage[all,2cell]{xy} \UseAllTwocells \SilentMatrices
\usepackage[T1]{fontenc}
\usepackage{xcolor}
\usepackage[toc,page]{appendix}
\usepackage[pdftex]{hyperref}
\usepackage{tikz}
\usepackage{tikz-cd}
\usepackage{blkarray}
\newtheorem{theo}{Theorem}[section]
\newtheorem{coro}[theo]{Corollary}
\newtheorem{prop}[theo]{Proposition}
\newtheorem{lemm}[theo]{Lemma}
\newtheorem{conj}[theo]{Conjecture}
\newtheorem{ques}[theo]{Question}

\theoremstyle{definition}
\newtheorem{defi}[theo]{Definition}

\theoremstyle{remark}
\newtheorem{rema}[theo]{Remark}

\newcommand{\lra}{\longrightarrow}

\DeclareMathOperator{\End}{End}
\DeclareMathOperator{\Hom}{Hom}

\DeclareMathOperator{\Ext}{Ext}

\newcommand{\id}{\mathrm{id}}

\newcommand{\HH}{\mathcal H}
\newcommand{\LL}{\mathcal L}
\newcommand{\PP}{\mathbb P}
\newcommand{\OO}{\mathcal O}

\newcommand{\vs}{\mathrm{vs}}

\newcommand{\dR}{\mathrm{dR}}
\newcommand{\Dol}{\mathrm{Dol}}
\newcommand{\Hod}{\mathrm{Hod}}

\DeclareMathOperator{\Conn}{Conn}
\DeclareMathOperator{\tr}{tr}

\DeclareMathOperator{\Pic}{Pic}

\title[Nonabelian Hodge theory for rank-2 bundles]{Nonabelian Hodge theory in positive characteristic for rank-2 vector bundles over curves}

\author{Christian Pauly}
\address{Laboratoire J.A. Dieudonnée, Université Côte d'Azur, Parc Valrose, 06108 Nice Cedex 02, France}
\email{christian.pauly@univ-cotedazur.fr}

\author{Thanh Tam Pham}
\address{Department of Mathematics, Hanoi Pedagogical University 2, Phu Tho, Vietnam}
\email{phamthanhtam@hpu2.edu.vn}



\date{}

\begin{document}

\begin{abstract}
We describe the Nonabelian Hodge Theory correspondence in positive characteristic between nilpotent Higgs bundles and flat bundles with nilpotent $p$-curvature, as constructed in the seminal paper by Ogus-Vologodsky, in the special case of 
rank-$2$ vector bundles over a smooth projective curve. In this set-up we can rephrase this correspondence as a shift of the extension class of a rank-$2$ vector bundle by a cohomology class encoding a $W_2$-lift of the curve. As an application, we obtain that a general rank-$2$ vector bundle of degree $0$ admits a connection with non-zero nilpotent $p$-curvature. 
\end{abstract}

\maketitle

\tableofcontents

\section{Introduction}

Let $X$ be a smooth projective and connected curve of genus $g \geq 2$ over an algebraically closed field $\Bbbk$. 
Very stable vector bundles were introduced by Drinfeld and Laumon \cite{La88} and have attracted considerable attention in the literature because of their importance in the study of the Hitchin systems (e.g. \cite{HH22}, \cite{DP09}). A
vector bundle $E$ over $X$ is called {\em very stable} if $E$ does not have nonzero nilpotent Higgs fields 
$\phi \in H^0(X, \mathrm{End}(E) \otimes \Omega_X)$ and {\em wobbly} otherwise. It is known that the locus of 
wobbly bundles $\mathcal{W}$ in the coarse moduli space $\mathcal{M}(X,r)$ of semi-stable rank-$r$ vector bundles
$$ \mathcal{W} = \{ [E] \in \mathcal{M}(X,r) \ | \ \exists  \phi\in H^0(X, \mathrm{End}(E) \otimes \Omega_X) \ \phi \ 
\text{nonzero and nilpotent} \}  $$
is of codimension $1$ (and conjectured to be of pure codimension $1$). It plays a significant role in the Donagi-Pantev approach to the geometric Langlands correspondence. 

\bigskip

In this paper we consider the case when the base field $\Bbbk$ has positive characteristic $p > 0$.
We denote the relative Frobenius by $F : X \longrightarrow X'$. Nonabelian Hodge theory in positive characteristic, as developped by Ogus-Vologodsky \cite{OV07}, Groechenig \cite{Gr16} and Chen-Zhu \cite{CZ15}, identifies the nilpotent cone consisting of pairs of nilpotent Higgs bundles $(E,\phi)$ over the Frobenius twist $X'$ 
with pairs $(V,\nabla)$ of flat bundles over $X$ having nilpotent $p$-curvature $\psi_\nabla$. This correspondence depends on 
the choice of a lifting of $X$ to the ring of Witt vectors $W_2(\Bbbk)$. Our motivating question is to 
investigate the different manifestations of ``wobbliness" for flat bundles $(E,\nabla)$. We restrict 
ourselves to rank-$2$ vector bundles, because in this case the geometry of the nilpotent cone of Higgs bundles is completely understood (see e.g \cite{PP21}). Some of our results are valid for higher rank bundles.

\subsection{The torsor of extension classes associated with $W_2(\Bbbk)$-liftings of the curve $X$} 
We fix a $W_2(\Bbbk)$-lifting of $X$ and choose 
local liftings of the relative Frobenius. Then the Deligne-Illusie Lemma \cite{DI87} asserts that the divided 
differences of two local liftings on an overlap form a \v{C}ech cocycle with values in $F^* T_{X'}$. We denote its class
by 
\begin{equation} \label{classh}
 h \in H^1(X, F^*T_{X'}).
\end{equation} 
We denote by $B_{X'}$ the sheaf of locally exact differential forms over $X'$ and consider the natural exact sequence
over $X'$ given by the Cartier operator
$$ 0 \lra B_{X'} \lra F_* \Omega_{X} \lra \Omega_{X'} \lra 0. $$  
This exact sequence determines an extension class 
$$[Ca] \in \mathrm{Ext}^1(\Omega_{X'}, B_{X'}) = H^1(X', B_{X'} \otimes T_{X'}).$$
Then we will show
\begin{lemm}
Varying the $W_2(\Bbbk)$-lifting of $X$ we obtain a set $\mathcal{L}(X)$ of 
classes $h$ as defined in \eqref{classh}.
Then $\mathcal{L}(X)$ coincides
with the fiber over the nonzero class $[Ca]$ of the linear map
$$ H^1(X, F^*T_{X'}) =  H^1(X', F_*\mathcal{O} \otimes T_{X'}) \lra H^1(X', B_{X'} \otimes T_{X'}).$$
\end{lemm}
Thus $\mathcal{L}(X)$ is an affine space modeled on the vector space $H^1(X',T_{X'})$ of first order infinitesimal
deformations of $X'$. We note that $h \in \mathcal{L}(X)$ is the gluing class entering the inverse Cartier transform
defined in \cite{OV07} and \cite{LSZ15}.

\subsection{The Dolbeault and de Rham nilpotent cones}
Let $\mathcal{M}(X',2)$ (resp. $\mathcal{M}(X,2)$) denote the coarse moduli space of semi-stable rank-$2$ vector bundles of trivial determinant over $X'$ (resp. $X$). These are irreducible projective varieties of dimension $3g-3$. 
We also consider $\mathcal{M}_{\Dol}(X',2)$ (resp. $\mathcal{M}_{\dR}(X,2)$) the coarse moduli spaces of semi-stable
rank-$2$ Higgs bundles (resp. flat bundles) with trivial determinant over $X'$ (resp. $X$). Then we have two Hitchin morphisms
\begin{equation}
\begin{tikzcd}
\mathcal{M}_{\Dol}(X',2) \ar[dr,"h_{\Dol}"] &  &  \ar[dl,swap,"h_{\dR}"] \mathcal{M}_{\dR}(X,2) \\
       & \mathcal{A}_{X'}^0 &  
\end{tikzcd}
\end{equation}
into the $\mathrm{SL}_2$-Hitchin base $\mathcal{A}_{X'}^0 := H^0(X', \Omega^{\otimes 2}_{X'})$. 
The first one $h_{\Dol}$ is the usual Hitchin morphism defined by $h_{\Dol}(E,\phi) = \det \phi$ and the second one $h_{\dR}$ is obtained by taking the Frobenius descent of 
$h_{\Dol}$ applied to the $p$-curvature $\psi_\nabla$, i.e. $F^*(h_{\dR}(V, \nabla)) = \det \psi_\nabla$. 
It is shown in \cite[Corollary 3.28]{Gr16} (see also \cite{dCGZ24}, \cite{CZ15}) that for every spectral data $a \in \mathcal{A}_{X'}^0$, a very good splitting of the Azumaya algebra of crystalline differential operators on the 
spectral curve $Y_a$ gives an isomorphism between the fibers
$$ h_{\Dol}^{-1}(a) \stackrel{\sim}{\lra} h_{\dR}^{-1}(a). $$
In this paper we will study in detail the isomorphisms obtained for $a=0$. Given an element $h \in \mathcal{L}(X)$ Ogus and Vologodsky
construct in \cite{OV07} an isomorphism, called the inverse Cartier transform, between the reduced underlying schemes 
\footnote{Strictly speaking \cite{OV07} constructs an equivalence of categories. Subsequently \cite{Gr16} and \cite{dCGZ24}
show that this equivalence of categories preserves semi-stability and induces an algebraic morphism}
$$ \mathrm{OV}_h : \mathcal{N}_{\Dol} := h_{\Dol}^{-1}(0)_{red}   \stackrel{\sim}{\lra}  \mathcal{N}_{\dR} := h_{\dR}^{-1}(0)_{red}.  $$
We will decompose (see Section 6) the Dolbeault and De Rham nilpotent cones $\mathcal{N}_{\Dol}$ and $\mathcal{N}_{\dR}$ into a union of
equidimensional components
$$ \mathcal{N}_{\Dol} = C'_0 \cup \bigcup_{k=2}^g C'_k \qquad \text{and} \qquad \mathcal{N}_{\dR} = C_0 \cup\bigcup_{k=2}^g C_k,$$
according to the degree of the kernel line subbundle of the nilpotent Higgs field (resp. $p$-curvature). Note that
for $k=0$ and $2 \leq k \leq g-1$ the components $C_k$ and $C'_k$ are also irreducible. We obtain therefore for each 
$k$ a restricted isomorphism $\mathrm{OV}_h : C'_k  \stackrel{\sim}{\lra} C_k $. 

The case $k= 0$ is easily described as follows. We have
$$ C'_0 = \mathcal{M}(X',2) \subset \mathcal{M}_{\Dol}(X',2)  \qquad \text{and} \qquad 
 C_0 = \{ (V, \nabla) \ | \ \psi_\nabla = 0 \} \subset \mathcal{M}_{\dR}(X,2), $$
 where $C'_0$ is the zero-Higgs-field component. Then $\mathrm{OV}_h (E,0) = (F^*E, \nabla_{can})$ and Cartier's theorem
 implies that $\mathrm{OV}_h$ is an isomorphism. Note that in this case the correspondence is independent of the
 class $h$, which stands in contrast with the other components.
 
\bigskip

Building on the work by Lan-Sheng-Zuo \cite{LSZ15}, we will give a description of  $\mathrm{OV}_h$ for the components 
$C'_k$ with $k \not= 0$ in terms of extension classes. It is shown in
\cite[Section 2]{PP21} that for $2 \leq k \leq g$ a general nilpotent Higgs bundle $(E, \phi) \in C'_k$ can be realized
as an extension 
$$ 0 \lra L \stackrel{i_E}{\lra} E \stackrel{\pi_E}{\lra} L^{-1} \lra 0, \qquad (e) $$
where $L$ is a line bundle of degree $\deg L = 1 -k$, the element $e \in H^1(L^2) = \mathrm{Ext}^1(L^{-1}, L)$ is its
extension class, and the Higgs field $\phi$ is obtained as composition $\phi = (i_E \otimes \mathrm{id}_{\Omega_{X'}}) \circ u \circ \pi_E$ for
a non-zero section $u \in H^0(L^2 \otimes \Omega_{X'}) = \mathrm{Hom}(L^{-1},L \otimes \Omega_{X'})$. We then show
\begin{prop} \label{prop:inverse-class}
The image $(V,\nabla) := \mathrm{OV}_h(E,\phi) $ is an extension
$$ 0 \lra F^*L \stackrel{i_V}{\lra} V \stackrel{\pi_V}{\lra} F^*L^{-1} \lra 0, \qquad (v) $$
with extension class $v \in H^1(F^* L^2)$ given by 
$$ v = F^*(e) + F^*(u) \smile h,$$
where $\smile$ denotes cup-product. The connection $\nabla$ on $V$ has nilpotent p-curvature 
$\psi_\nabla = - i_V \circ F^*(u) \circ \pi_V \in H^0 (X, \mathrm{End}_0(V) \otimes F^*\Omega_{X'})$.
\end{prop}

\subsection{The vector space of $t$-connections on a fixed bundle $V$ and the restricted de Rham--Hitchin morphism}

We now consider for a stable rank-$2$ vector bundle $V$ over $X$ the affine space of connections $\mathrm{Conn}^0(V)$ on $V$ with fixed trivial determinant, i.e. $\det \nabla = d$, and the restricted de Rham--Hitchin morphism
$$h_{dR,V} : \mathrm{Conn}^0(V) \lra \mathcal{A}^0_{X'},$$
which, by definition, is the restriction of $h_{\dR}$ to $\mathrm{Conn}^0(V) \subset \mathcal{M}_{\dR}(X,2)$.
Clearly, as a set, the fiber $h_{dR,V}^{-1}(0)$ identifies with the set of connections $\nabla$ on $V$ having nilpotent 
$p$-curvature. We then show 
\begin{prop} \label{prop:stabledRfinite}
If $V$ is very stable, then the morphism $h_{dR,V}$ is finite. 
\end{prop}

We don't know whether the converse holds, as for Higgs bundles --- recall 
\cite{PP18} that for the restricted Dolbeault--Hitchin morphism $h_{Dol,V}$ we have an equivalence between $V$ very stable
and $h_{Dol,V}$ finite. The strategy of the proof of Proposition \ref{prop:stabledRfinite} is to consider the vector space
of all $t$-connections $t\mathrm{Conn}^0(V)$ on $V$ and to apply a result by H. Zelaci \cite{Ze18} to the restricted Hodge--Hitchin morphism from $t\mathrm{Conn}^0(V)$ to $\mathcal{A}^0_{X'} \times \Bbbk$ given by homogeneous
polynomials.

\subsection{Dominance of the de Rham forgetful map}

Finally, we consider the rational de Rham forgetful map 
$$ \pi_{\dR}: \mathcal{N}_{\dR}^{\not= 0} = \bigcup_{k=2}^g C_k  \dashrightarrow \mathcal{M}(X,2), \qquad (V,\nabla) \mapsto V.$$
Since $\dim C_k = \dim \mathcal{M}(X,2)$ for every $k$, it is natural to ask whether $\pi_{\dR}$ is dominant. 
\begin{prop} \label{dominance}
There exists an integer $k$ satisfying  $2 \leq k \leq g$ such that the 
restriction of $\pi_{\dR}$ to the component $C_k$ is dominant, hence generically finite.
\end{prop}

The proof of Proposition \ref{dominance} relies on the computation of the length 
of the $0$-dimensional scheme $h^{-1}_{dR,V}(0)$ for a 
general bundle $V$. We then exclude the case that the reduced scheme $h^{-1}_{dR,V}(0)_{red}$ is entirely contained in the zero-$p$-curvature component $C_0$ using an estimate of the degree of the generalized Verschiebung.

\bigskip

As an immediate corollary of Proposition \ref{dominance} we obtain

\begin{coro}
Over any curve $X$ a general rank-$2$ vector bundle of degree $0$ admits a connection $\nabla$ such that its 
$p$-curvature $\psi_\nabla$ is nonzero and nilpotent.
\end{coro}

Unfortunately, we are not able to determine the subset of integers $k$ for which the restriction of 
$\pi_{\dR}$ to $C_k$ is dominant. We note that our partial result contrasts with its analogue result 
for Higgs bundles: the rational Dolbeault forgetful map
$$ \pi_{\Dol} :  \mathcal{N}_{\Dol}^{\not= 0} = \bigcup_{k=2}^g C'_k  \dashrightarrow \mathcal{M}(X',2), \qquad (E,\phi) \mapsto E$$
is 
completely described in \cite{PP21} and we know that for $\max(2,\frac{g}{2}) \leq k \leq g$ the images 
$\pi_{\Dol}(C'_k)$ are divisors, called wobbly divisors, and for $2 \leq k < \frac{g}{2}$ the images
$\pi_{\Dol}(C'_k)$ have codimension at least $2$.

For simplicity of presentation we state and prove our main results only for rank-$2$ vector bundles with fixed trivial determinant. We will address possible generalization to higher rank vector bundles in a future paper.  
At the end of this paper (Section 7) we add a list of open questions.

\bigskip
\textbf{Acknowledgements.} The first author wishes to thank the organizers of the conference "Arithmetic -- Algebraic Geometry" held in Nha Trang (Vietnam) in June 2025, where this project has been started. The second author thanks the IRL France-Vietnam and the University of Montpellier for financial support for a research visit in Montpellier (France) in October 2026.  We also thank Q.D. Dinh, H.H. Phung, M. Sheng, Y. Wakabayashi, S. Zhang and  K. Zuo for useful discussions. The second author's research is funded by NAFOSTED under grant number 101.04-2025.20 “Some problems on fundamental groups, cohomology, and moduli space”.

\section{Preliminaries on $W_2(\Bbbk)$-liftings, Frobenius lifts and Deligne-Illusie Lemma}

%In this section, denote by \(\mathcal M(X',2,d)\) the moduli space of semistable rank-$2$ bundles of degree \(d\) on \(X'\), and \(\mathcal M(X,2,pd)\) %the corresponding moduli space of degree \(pd\) on \(X\). 

\subsection{Local Frobenius lifts and the Deligne--Illusie lemma}

Let $\Bbbk$ be a perfect field of characteristic $p>0$, and let $F:X\to X'$ be the relative Frobenius of 
the curve $X$. Fix a $W_2(\Bbbk)$-lifting $\widetilde X$ of $X$, and denote its Frobenius twist by $\widetilde X'$. We recall the local construction of \cite[proof of Theorem~2.1, (b)--(d)]{DI87}; see also \cite[Section~2.1]{LSZ15}.

Choose an affine open cover $X=\bigcup_\alpha U_\alpha$ and local $W_2(\Bbbk)$-morphisms $\widetilde F_\alpha\colon \widetilde U_\alpha\longrightarrow\widetilde U'_\alpha$ lifting $F|_{U_\alpha}$. Their differentials are divisible by $p$ and define $\OO_{U_\alpha}$-linear maps
\[
\zeta_\alpha = \frac{d\widetilde F_\alpha}{p}\bmod p
\colon
F^*\Omega^1_{U'_\alpha} \longrightarrow\Omega^1_{U_\alpha}.
\]
On overlaps $U_{\alpha \beta} = U_\alpha \cap U_\beta$, the divided differences $\dfrac{\widetilde F_\beta^\#- \widetilde F_\alpha^\#}{p}\bmod p$ are Frobenius derivations and hence induce $\OO_{U_{\alpha\beta}}$-linear maps
\[
h_{\alpha\beta}\colon F^*\Omega^1_{U'_{\alpha\beta}} \longrightarrow\OO_{U_{\alpha\beta}}.
\]
Here division by $p$ uses the canonical identification $p\widetilde M\simeq\widetilde M/p\widetilde M$ for flat $W_2(\Bbbk)$-modules $\widetilde M$.

\begin{lemm}[Deligne--Illusie]\label{lemma-deligneIlucsse}
	The preceding maps satisfy
	\[
	\zeta_\beta-\zeta_\alpha=dh_{\alpha\beta}, \qquad h_{\alpha\beta}+h_{\beta\gamma}=h_{\alpha\gamma}.
	\]
	The first identity is understood on $F^{-1}\Omega^1_{U'_{\alpha\beta}}$, and the second holds on triple overlaps.
\end{lemm}

%%%%%%%%%%%%%%%%%%%%%%%%%%%%%%%%%%%%%%%%%%%%%%%%%%%%%%%%%%%%%%%
\subsection{Extension classes and transition functions} Let $E',E''$ be vector bundles on $X$ and
\begin{equation}\label{extension-01}
0 \longrightarrow E' \xrightarrow{\;i\;} E \xrightarrow{\;\pi\;} E'' \longrightarrow 0
\end{equation}
be an extension of $E''$ by $E'$. Following \cite[Section 1, p.185]{At57}, there exists an affine covering $\big\{U_\alpha\big\}$ of $X$ over which the extension is trivial, i.e., the exact sequence \[0 \longrightarrow E'_\alpha \xrightarrow{\;i\;} E_\alpha \xrightarrow{\;\pi_\alpha\;} E''_\alpha \longrightarrow 0\]
admits a section $E''_\alpha \xrightarrow{s_\alpha} E_\alpha$ of $\pi_\alpha$, where we write $E_\alpha = E|_{U_\alpha}$ (same for $E'_\alpha, E''_\alpha$).
Then, over $U_{\alpha\beta} = U_\alpha \cap U_\beta$, we have $s_{\alpha\beta} := s_\beta - s_\alpha \in \operatorname{Hom}\!\left(E''_{\alpha\beta}, E'_{\alpha\beta}\right)$. Therefore, $\{ s_{\alpha\beta} \}$ is a $1$-cocycle with values in $\operatorname{Hom}(E'', E')$ which defines the extension class
\[
[E]=[s_{\alpha\beta}] \in H^1\!\left(\operatorname{Hom}(E'', E')\right).
\]

Following \cite[Section 1, p. 185]{At57}, we consider the isomorphism $\varphi_\alpha \colon E'_\alpha \oplus E''_\alpha \longrightarrow E_\alpha$ over $U_\alpha$ given by $(e', e'') \longmapsto e' + s_\alpha(e'')$. Over $U_{\alpha\beta}$, we define the transition functions $g_{\alpha\beta} := \varphi_\alpha^{-1} \circ \varphi_\beta$ by
\[
g_{\alpha\beta} \colon E'_\beta \oplus E''_\beta \xrightarrow{\ \varphi_\beta\ } E_{\alpha\beta} \xrightarrow{\ \varphi_\alpha^{-1}\ } E'_\alpha \oplus E''_\alpha,
\]
where $\varphi_\beta(e', e'') = e' + s_\beta(e'')$ and $\varphi_\alpha^{-1}\big(e' + s_\beta(e'')\big) = \bigl(e' + s_\beta(e'') - s_\alpha(e''),\, e''\bigr)$.
Explicitly, we have
\begin{align*}
    g_{\alpha\beta}(e', e'')= \bigl(e' + (s_\beta - s_\alpha)(e''),\, e''\bigr)= \bigl(e' + s_{\alpha\beta}(e''),\, e''\bigr).
\end{align*}
Then, the vector bundle $E$ is
\[
E = \left( \bigsqcup_{\alpha} \, E'_\alpha \oplus E''_\alpha \right) \big/ (g_{\alpha\beta}).
\]

\begin{rema}
Let $G_{\alpha\beta} \colon E_{\alpha\beta} \longrightarrow E_{\alpha\beta}$ be transition functions, assume that
\[
G_{\alpha\beta} = \id_{E_{\alpha\beta}} + i \circ u_{\alpha\beta} \circ \pi,
\]
where $u_{\alpha\beta} \in \operatorname{Hom}\!\left(E''_{\alpha\beta}, E'_{\alpha\beta}\right)$ is a \(1\)-cocycle. According to \cite[Proposition 2]{At57}, the class $[u_{\alpha\beta}] \in H^1\!\left(X, \operatorname{Hom}\bigr(E'', E'\bigr)\right)$ defines a extension class $[u_{\alpha\beta}]\in \operatorname{Ext}^1\bigr(E'', E'\bigr)$. Notice that the cocycle condition $u_{\alpha\beta} + u_{\beta\gamma} = u_{\alpha\gamma}$ is equivalent to
\[
G_{\alpha\beta} \circ G_{\beta\gamma} = G_{\alpha\gamma}.
\]

Let $V$ be the vector bundle obtained from $E$ by gluing along $G_{\alpha\beta}$. Then, $V$ is by definition
\[
V = \bigsqcup_{\alpha} E_{\alpha} \big/ (G_{\alpha\beta}),
\]
where $G_{\alpha\beta} : E_\beta\big|_{U_{\alpha\beta}} \longrightarrow E_\alpha\big|_{U_{\alpha\beta}}$ given by $(x,e) \longmapsto \bigl(x,\, G_{\alpha\beta}(x)(e)\bigr)$ for $x \in U_{\alpha\beta}$.
\end{rema}

\begin{prop}\label{extension-sum of classes}
Let $0 \longrightarrow E' \xrightarrow{\;i\;} E \xrightarrow{\;\pi\;} E'' \longrightarrow 0$ be an extension of $E''$ by $E'$ over $X$, defined with respect to an open cover $\big\{U_\alpha\big\}$ of $X$ by the transition functions
$g_{\alpha \beta}$. Let $G_{\alpha\beta} \colon E_{\alpha\beta} \longrightarrow E_{\alpha\beta}$ be transition functions of the form$$G_{\alpha\beta} = \mathrm{id}_{E_{\alpha\beta}} + i \circ u_{\alpha\beta} \circ \pi,$$where $u_{\alpha\beta} \in \operatorname{Hom}\!\left(E''_{\alpha\beta}, E'_{\alpha\beta}\right)$ is a $1$-cocycle, and let $V$ be the vector bundle obtained from $E$ by gluing via $G_{\alpha\beta}$. Then
\[
[V]=[E]+[u_{\alpha\beta}]\in \operatorname{Ext}^1\bigr(E'', E'\bigr).
\]
\end{prop}
\begin{proof}
The transition functions $g'_{\alpha\beta}=\varphi_\alpha^{-1} \circ G_{\alpha\beta} \circ \varphi_\beta$ over $U_{\alpha\beta}$ are given by the composition
\[
g'_{\alpha\beta}:
E'_\beta \oplus E''_\beta
\xrightarrow{\;\varphi_\beta\;} E_\beta \xrightarrow{\;G_{\alpha\beta}\;} E_\alpha
\xrightarrow{\;\varphi_\alpha^{-1}\;} E'_\alpha \oplus E''_\alpha .
\]
A direct computation yields
\begin{align*}
\varphi_\beta(e',e'') &= e' + s_\beta(e''),\\
G_{\alpha\beta}\bigl(e' + s_\beta(e'')\bigr) &= e' + s_\beta(e'') + u_{\alpha\beta}(e''),\\
g'_{\alpha \beta}(e',e'') =\varphi_\alpha^{-1}\bigl(e' + s_\beta(e'') + u_{\alpha\beta}(e'')\bigr) &= \bigl(e' + s_\beta(e'') - s_\alpha(e'') + u_{\alpha\beta}(e''), e''\bigr),
\end{align*}
whence we obtain the following commutative diagram over $U_{\alpha\beta}$:
$$
\begin{tikzcd}
E'_\beta \oplus E''_\beta \arrow[r, "\varphi_\beta"] \arrow[d, "g'_{\alpha\beta}"'] & E_\beta \arrow[d, "G_{\alpha\beta}"]\\
E'_\alpha \oplus E''_\alpha \arrow[r, "\varphi_\alpha"'] & E_\alpha
\end{tikzcd}
$$

Therefore, the transition matrix of $V$ is
\[
g'_{\alpha\beta}
    =\varphi_\alpha^{-1} \circ G_{\alpha\beta} \circ \varphi_\beta\\
    =\begin{pmatrix}
        \id& s_{\alpha\beta}+u_{\alpha\beta}\\
        0&\id
    \end{pmatrix}.
\]
This shows that $s_{\alpha\beta}+u_{\alpha\beta}$ is the $1$-cocycle defined by $V$, hence
    \[
    [V]=[s_{\alpha\beta}+u_{\alpha\beta}]=[s_{\alpha\beta}]+[u_{\alpha\beta}]=[E]+[u_{\alpha\beta}].    
    \]
\end{proof}

%%%%%%%%%%%%%%%%%%%%%%%%%%%%%%%%%%%%%%%%%%%%%%%%%%%%%%%%%%%%%%%
\subsection{Relation to the sheaf of locally exact differentials}\label{subsec:Relation to exact differentials}

The Deligne--Illusie cocycle admits a useful interpretation in terms of the Cartier sequence. Fix a $W_2(\Bbbk)$-lifting of $X$ and choose local Frobenius liftings on an open cover $\{U_\alpha\}$. Their differences determine a cocycle
\[
h_{\alpha\beta}\colon
F^*\Omega_{U'_{\alpha\beta}} \longrightarrow \OO_{U_{\alpha\beta}}.
\]
Since $F$ is finite, adjunction and the projection formula give canonical identifications
\begin{align*}
	H^1\bigl(X,(F^*\Omega_{X'})^\vee\bigr) &\simeq
	H^1\bigl(X',F_*((F^*\Omega_{X'})^\vee)\bigr)\simeq
	H^1\bigl(X',(F_*\OO_X)\otimes T_{X'}\bigr),
\end{align*}
where we denote $T_{X'} = \Omega_{X'}^\vee$.
We write $h=[h_{\alpha\beta}]$ for the resulting class, viewing each $h_{\alpha\beta}$ also as its adjoint $$\Omega_{U'_{\alpha\beta}} \to F_*\OO_{U'_{\alpha\beta}}.$$

Let $B_{X'}:=\operatorname{im} \bigl(d\colon F_*\OO_X\longrightarrow F_*\Omega_{X}\bigr)$ be the sheaf of locally exact differentials on $X'$. As $X$ is a smooth curve, the Cartier operator yields the exact sequences
\begin{equation}\label{ex:exact-functions}
	0\longrightarrow\OO_{X'} \longrightarrow F_*\OO_X \xrightarrow{d} B_{X'} \longrightarrow0
\end{equation}
and
\begin{equation}\label{s1:cartier}
	0\longrightarrow B_{X'} \longrightarrow F_*\Omega_{X} \xrightarrow{C}\Omega_{X'} \longrightarrow0.
\end{equation}
Denote the extension class of \eqref{s1:cartier} by
\[
[Ca]\in
\operatorname{Ext}^1_{X'}\left(\Omega_{X'},B_{X'}\right) \simeq H^1\bigl(X',B_{X'} \otimes T_{X'} \bigr).
\]
Tensoring \eqref{ex:exact-functions} with $T_{X'}$ and using $H^2(X',T_{X'})=0$, since $\dim X'=1$, 
gives a surjection
\[
\pi\colon
H^1\bigl(X',(F_*\OO_X)\otimes T_{X'}\bigr) \longrightarrow H^1\bigl(X',B_{X'} \otimes T_{X'} \bigr).
\]

To identify $\pi(h)$, let $\zeta_\alpha\colon \Omega_{U'_\alpha}\to F_*\Omega_{U'_\alpha}$ 
be the local Cartier splittings induced by the chosen Frobenius liftings. With the difference convention $d h_{\alpha\beta} =\zeta_\beta-\zeta_\alpha$, Lemma \ref{lemma-deligneIlucsse}  and the definition of the extension class (Section 2.2) shows that $\{d h_{\alpha\beta}\}$ is the cocycle representing the extension \eqref{s1:cartier}. Consequently, we obtain the relation
$$\pi(h) = [Ca],$$
which provides a cohomological formulation of the compatibility between the Deligne--Illusie class and the Cartier sequence. In fact, the surjectivity of $\pi$ guaranties the existence of a class $h \in H^1\bigl(X', (F_*\OO_X)\otimes T_{X'} \bigr)$ that satisfies this identity independently of any $W_2(\Bbbk)$-lifting of $X$. Therefore, for arguments relying solely on this compatibility, one may work directly with such a class without specifying a 
$W_2(\Bbbk)$-lifting from which it arises.

%%%%%%%%%%%%%%%%%%%%%%%%%%%%%%%%%%%%%%%%%%%%%%%%%%%%%%%%%%%%%%%
\section{The Ogus-Vologodsky correspondence for rank-2 bundles via extension classes}
\subsection{The affine space of lifts of the Cartier class}

We now study the fibre of $\pi$ introduced in Subsection \ref{subsec:Relation to exact differentials} over the Cartier class $[Ca]$. Via the canonical identification
\begin{equation}\label{eq:adjunction-cohomology}
	H^1(X,F^*T_{X'}) \simeq H^1(X',F_*\OO_X\otimes T_{X'}),
\end{equation}
we regard this fibre as a subset of $H^1(X,F^*T_{X'})$.

\begin{defi}\label{def:lift-space}
We define the cohomological lift space of the Cartier class to be the affine space
\[
\mathcal L(X):=\pi^{-1}([Ca]) \subset H^1(X,F^*T_{X'}).
\]
We also define the vector space 
$$ \widehat{\mathcal{L}(X)} := \pi^{-1}(\Bbbk[Ca]) \subset H^1(X,F^*T_{X'}).$$
\end{defi}

\begin{lemm}\label{lem:Fro-liftings}
	The set $\mathcal{L}(X)$ is nonempty and naturally forms an affine space modeled on $H^1(X', T_{X'})$ (or more precisely, a torsor under $H^1(X', T_{X'})$), with the translation action given by
	\[
	H^1(X',T_{X'})\times\mathcal L(X) \longrightarrow\mathcal L(X), \qquad (\eta,h)\longmapsto h+F^*\eta.
	\]
	In particular,
	\[
	\dim\mathcal L(X)=3g-3.
	\]
\end{lemm}

\begin{proof}
 Inclusion $B_{X'}\hookrightarrow F_*\Omega_{X}$ and projection formula give an injection
	\[
	H^0(X',B_{X'}\otimes T_{X'}) \hookrightarrow H^0(X,\Omega_{X}\otimes F^*T_{X'}).
	\]
As $\deg(\Omega_{X}\otimes F^*T_{X'}) = (1-p)(2g-2) < 0$, both spaces vanish. Moreover, since $H^2(X',T_{X'})=0$, tensoring \eqref{ex:exact-functions} with $T_{X'}$ and taking cohomology yields the short exact sequence
	\begin{equation}\label{eq:cohomology-sequence}
		0\longrightarrow H^1(X',T_{X'})
		\xrightarrow{F^*}H^1(X,F^*T_{X'})
		\xrightarrow{\pi}H^1(X',B_{X'}\otimes T_{X'})
		\longrightarrow0.
	\end{equation}
	Thus, $\mathcal L(X)$ is a nonempty fibre of $\pi$, and $\ker\pi=F^*H^1(X',T_{X'})$ acts simply transitively on it. The dimension formula follows from Riemann--Roch and 	$H^0(X',T_{X'})=0$.
\end{proof}

As explained in Subsection \ref{subsec:Relation to exact differentials}, the cocycle associated with the Frobenius-lifting data determines a class $h=[h_{\alpha\beta}]\in\mathcal L(X)$. In the constructions below, we work directly with a point of $\mathcal L(X)$, and it is unnecessary to choose a $W_2(\Bbbk)$-lifting that realizes that point.

%%%%%%%%%%%%%%%%%%%%%%%%%%%%%%%%%%%%%%%%%%%%%%%%%%%%%%%%%%%%%%%
\subsection{The normalized connection on \texorpdfstring{$\HH_h^2$}{H2}}

Let $h\in\mathcal L(X)$, where $\mathcal L(X)$ is the affine space defined in Definition~\ref{def:lift-space}. We realize $h$ as an extension class and use its defining property $\pi(h)=[Ca]$ to construct the normalized connection on $\HH_h^2$.

%\paragraph{The two extension descriptions.}
Since $\Omega_{X'}$ is locally free, identification \eqref{eq:adjunction-cohomology} induces a canonical isomorphism
\begin{equation}\label{eq:extension-adjunction}
	\operatorname{Ext}^1_X(F^*\Omega_{X'},\OO_X)\simeq \operatorname{Ext}^1_{X'}(\Omega_{X'},F_*\OO_X).
\end{equation}
Accordingly, $h$ determines two corresponding extensions (with the superscripts recording the ranks)
\begin{equation} 
    0\longrightarrow\OO_X \xrightarrow{i_h}\HH_h^2 \xrightarrow{\pi_h}F^*\Omega_{X'} \longrightarrow0, \label{eq:H2}
\end{equation}
\begin{equation}
    0\longrightarrow F_*\OO_X \xrightarrow{i'_h}\HH_h^{p+1} \xrightarrow{\pi'_h}\Omega_{X'} \longrightarrow0. \label{eq:Hp1}
\end{equation}

Let $\eta\colon\Omega_{X'}\longrightarrow F_*F^*\Omega_{X'}$ be defined by the adjunction of the identity and $\operatorname{ev}_h\colon F^*F_*\OO_X\longrightarrow\OO_X$ be the evaluation map. Since $F_*F^*\Omega_{X'}=F_*\OO\otimes\Omega$, by applying $F_*$ to \eqref{eq:H2} and pulling back the resulting extension along $\eta$, we obtain
\[
\HH_h^{p+1}\simeq F_*\HH_h^2 \times_{F_*F^*\Omega_{X'}}\Omega_{X'}.
\]
Conversely, pulling \eqref{eq:Hp1} back by $F$ and pushing out along $\operatorname{ev}_h$ gives
\[
\HH_h^2\simeq F^*\HH_h^{p+1} \sqcup_{F^*F_*\OO_X}\OO_X.
\]
Here, exactness is preserved because $F$ is finite and flat. These two constructions are inverse to each other and yield natural maps
\[
\HH_h^{p+1}\hookrightarrow F_*\HH_h^2,\qquad F^*\HH_h^{p+1}\twoheadrightarrow\HH_h^2.
\]

\begin{lemm}\label{lem:normalized-connection}
	For every $h\in\mathcal L(X)$, the extension \eqref{eq:H2} admits a unique connection $\nabla$ making $i_h$ and $\pi_h$ horizontal for the canonical connections on $\OO_X$ and $F^*\Omega_{X'}$. It is integrable, and its $p$-curvature is
	\begin{equation}\label{eq:normalized-p-curvature}
		\psi_{\nabla}	=-(i_h\otimes\id_{F^*\Omega_{X'}})\circ \pi_h.
	\end{equation}
	In particular, its $p$-curvature is nonzero and nilpotent of order $2$.
	%\begin{equation}\label{eq:p-curvature-convention}		\psi_\nabla(D)=\nabla_D^p-\nabla_{D^{[p]}}.	\end{equation}
\end{lemm}

\begin{proof}
	Under \eqref{eq:extension-adjunction}, the map $\pi$ is induced by pushing out extensions along $d\colon F_*\OO_X\to B_{X'}$. Thus, $\pi(h)=[Ca]$ says precisely that the pushout of \eqref{eq:Hp1} along $d$ is isomorphic, as an extension, to the Cartier sequence \eqref{s1:cartier}. Consequently, there is a commutative diagram
	\begin{equation}\label{eq:normalized-Cartier-diagram}
		\begin{tikzcd}[column sep=small]
			0 \arrow[r] & F_*\OO_X \arrow[r,"i'_h"] \arrow[d,"d"] & \HH_h^{p+1} \arrow[r,"\pi'_h"] \arrow[d,"\rho_h"] & \Omega_{X'} \arrow[r] \arrow[d,equal] &
			0 \\
			0 \arrow[r] & B_{X'} \arrow[r] & F_*\Omega_{X} \arrow[r,"C"] & \Omega_{X'} \arrow[r] & 0.
		\end{tikzcd}
	\end{equation}
	The map $\rho_h$ is unique, since the difference of two such maps factors through an element of 
	\[
	\operatorname{Hom}_{X'}(\Omega_{X'},B_{X'}) = H^0(X',B_{X'}\otimes T_{X'})=0.
	\]
	Note that this vanishing is established in the proof of Lemma~\ref{lem:Fro-liftings}.
	
	Choose an affine cover $\{U'_\alpha\}$ and local splittings $r_\alpha$ of $\pi'_h$. Write $r_\beta-r_\alpha=i'_h\bar h_{\alpha\beta}$ and $\bar h_{\alpha\beta}\colon\Omega_{X'}\longrightarrow F_*\OO_X$. Then, their adjoints $h_{\alpha\beta}\colon F^*\Omega_{X'}\to\OO_X$ represent $h$. Set
	\[
	\zeta_{h,\alpha}:=\rho_h r_\alpha \colon\Omega_{U'_\alpha} \longrightarrow F_*\Omega_{U_\alpha}.
	\]
	By \eqref{eq:normalized-Cartier-diagram},
	\begin{equation}\label{eq:adapted-Cartier-splittings}
		C\zeta_{h,\alpha}=\mathrm{id},\qquad \zeta_{h,\beta}-\zeta_{h,\alpha}=d\circ\bar h_{\alpha\beta}.
	\end{equation}
	In particular, these local Cartier splittings are chosen to be compatible with the representative of $h$. An arbitrary initial choice of Cartier splittings need not satisfy the second identity for that representative.
	
	Put $U_\alpha=F^{-1}(U'_\alpha)$. The corresponding local splittings of \eqref{eq:H2} identify $\HH_h^2|_{U_\alpha}$ with $(\OO_X\oplus F^*\Omega_{X'})|_{U_\alpha}$. With the convention that transition maps send coordinates on $U_\beta$ to coordinates on $U_\alpha$, the gluing maps are
	\[
	g_{\alpha\beta} =
    \begin{pmatrix}
		1&h_{\alpha\beta}\\
		0&1
	\end{pmatrix}.
	\]
	Denote the adjoint of $\zeta_{h,\alpha}$ by $\zeta_{h,\alpha}^{\sharp}\colon F^*\Omega_{U'_\alpha} \longrightarrow\Omega_{U_\alpha}$. Define
	\begin{equation}\label{eq:H-connection}
		\nabla_{h,\alpha} = 
        \begin{pmatrix}
			d&\zeta_{h,\alpha}^{\sharp}\\
			0&\nabla^{\mathrm{can}}
		\end{pmatrix}.
	\end{equation}
	In local horizontal frames for the canonical connections, write $\nabla_{h,\alpha}=d+A_\alpha$, where
	\[
	A_\alpha= 
    \begin{pmatrix}
		0&\zeta_{h,\alpha}^{\sharp}\\
		0&0
	\end{pmatrix}.
	\]
	Since \eqref{eq:adapted-Cartier-splittings}, $\zeta_{h,\beta}^{\sharp}-\zeta_{h,\alpha}^{\sharp}=dh_{\alpha\beta}$. Consequently,
	\begin{align*}
		g_{\alpha\beta}^{-1}A_\alpha g_{\alpha\beta}+g_{\alpha\beta}^{-1}dg_{\alpha\beta}&=
		\begin{pmatrix}
			0&\zeta_{h,\alpha}^{\sharp}\\
			0&0
		\end{pmatrix}
		+
		\begin{pmatrix}
			1&-h_{\alpha\beta}\\
			0&1
		\end{pmatrix}
		\begin{pmatrix}
			0&dh_{\alpha\beta}\\
			0&0
		\end{pmatrix}=A_\beta.
	\end{align*}
Thus, the local connections glue to a connection $\nabla$ on $\HH_h^2$. Both end terms of \eqref{eq:H2} have canonical connections, and integrability follows from $\bigwedge^2\Omega_{X}=0$.
	
	For completeness, we compute the sign in \eqref{eq:normalized-p-curvature}; compare the auxiliary connection in \cite[Section~2.3]{LSZ15}. Choose an etale coordinate $z$, with the corresponding coordinate $z'$ on $X'$, and put $D=\partial/\partial z$. After refining the cover, Cartier normalization gives
	\[
	\zeta_{h,\alpha}(dz')=z^{p-1}dz+db_\alpha.
	\]
	Let $e_1$ be the unit section of the subbundle, and let $e_2$ be the lift, in the chosen local splitting, of $F^*dz'$. Then
	\[
	\nabla_{h,D}e_1=0,\qquad \nabla_{h,D}e_2=(z^{p-1}+Db_\alpha)e_1.
	\]
	Since $D^{[p]}=D^p=0$,
	\[
	\psi_{\nabla}(D)e_2=D^{p-1}(z^{p-1}+Db_\alpha)e_1=-e_1.
	\]
	This proves \eqref{eq:normalized-p-curvature}. The image lies in the subbundle, which the $p$-curvature annihilates, so its values have pairwise zero products.
	
	Finally, the difference between two connections making $i_h$ and $\pi_h$ horizontal factors through a homomorphism $F^*\Omega_{X'}\to\Omega_{X}$. The space of such homomorphisms is
	\[
	H^0(X,\Omega_{X}\otimes F^*T_{X'})=0,
	\]
	because the line bundle has degree $(1-p)(2g-2)<0$. This proves uniqueness and independence of local choices.
\end{proof}

We call $\nabla$ the \textit{normalized connection} on $\HH_h^2$. Its construction uses $\pi(h)=[Ca]$ directly and does not require the choice of a $W_2(\Bbbk)$-lifting realizing $h$.

\begin{rema}\label{rem:normalization}
	Let $u\in H^1(X, F^*T_{X'})$ be arbitrary and suppose that its rank-$2$ extension carries a connection inducing the canonical connections on its subbundle and quotient. Write its local off-diagonal terms, by adjunction, as maps $\varphi_\alpha\colon\Omega_{X'}\to F_*\Omega_{X}$. The gluing condition gives $\varphi_\beta-\varphi_\alpha=d\circ\bar u_{\alpha\beta}$. Applying $C$ shows that the maps $C\varphi_\alpha$ glue to an endomorphism of $\Omega_{X'}$. Hence, $C\varphi_\alpha=\lambda\,\mathrm{id}$ for a scalar
	$\lambda\in H^0(X',\OO_{X'})=\Bbbk$, and 
    \begin{equation}\label{eq:scaled-Cartier-class}
		\pi(u)=\lambda[Ca].
	\end{equation}
	Conversely, \eqref{eq:scaled-Cartier-class} is sufficient for the existence of such a connection, by the same pushout argument.
	
	If $\lambda\ne0$, simultaneous rescaling of the cocycle and the off-diagonal terms gives the normalized construction for the class $\lambda^{-1}u\in\mathcal L(X)$. Rescaling the off-diagonal terms alone does not preserve the gluing relation for the original cocycle.
	
	If $\lambda=0$, then $\pi(u)=0$. By using \eqref{eq:cohomology-sequence}, there exist a unique $\eta\in H^1(X',T_{X'})$ such that $u=F^*\eta$. Thus, $\HH_u^2$ is the Frobenius pullback, as an extension, of an extension of $\Omega_{X'}$ by $\OO_{X'}$. The connection compatible with the canonical connections on its end terms is its canonical pullback connection. This uniqueness concerns connections compatible with the specified extension structure.
\end{rema}

%%%%%%%%%%%%%%%%%%%%%%%%%%%%%%%%%%%%%%%%%%%%%%%%%%%%%%%%%
\subsection{The inverse Cartier transform}

Fix $h\in\LL(X)$, together with normalized local data $(h_{\alpha\beta},\zeta_\alpha)$ as in the preceding subsections. In particular,
\[
dh_{\alpha\beta}=\zeta_\beta-\zeta_\alpha, \quad C\circ\zeta_\alpha=\id_{\Omega_{X'}}.
\]
Here, as before, we identify $\zeta_\alpha$ with its adjoint map $F^*\Omega_{X'}\to\Omega_X$ whenever convenient.

Let $(E,\phi)\in\mathcal M_{\Dol}(X',2)$ be a rank-$2$ Higgs bundle with trivial determinant, with $\phi\ne0$ and $\phi^2=0$. Put $L=\ker\phi$. The nilpotence of $\phi$ implies that its image lies in $L\otimes\Omega_{X'}$. Using the identification $E/L\simeq L^{-1}$, we therefore obtain \begin{equation}\label{eq:Higgs-extension}
	0\longrightarrow L\xrightarrow{i_E}E 	\xrightarrow{\pi_E}L^{-1}\longrightarrow0, \qquad \phi=(i_E\otimes\id_{\Omega_{X'}})\circ u\circ \pi_E,
\end{equation}
where
\[
0\ne u\in H^0(X',L^2\otimes\Omega_{X'}), \qquad [E]\in H^1(X',L^2).
\]
Throughout this subsection, $[E]$ denotes the class of the extension \eqref{eq:Higgs-extension}.

Pulling back \eqref{eq:Higgs-extension} by $F$ gives an extension with class $F^*[E]\in H^1(X,F^*L^2)$. On the other hand, $h\in H^1\bigl(X,(F^*\Omega_{X'})^\vee\bigr)$ and $F^*u\in H^0\bigl(X,F^*L^2\otimes F^*\Omega_{X'}\bigr)$. The cup product followed by contraction defines
\[
H^1\bigl(X,(F^*\Omega_{X'})^\vee\bigr)\otimes H^0\bigl(X,F^*L^2\otimes F^*\Omega_{X'}\bigr) \longrightarrow H^1(X,F^*L^2).
\]
We denote the resulting class by
\begin{equation}\label{eq:contraction}
	h\smile F^*u\in H^1(X,F^*L^2).
\end{equation}
More explicitly, it is represented by the cocycle
\[
\nu^h_{\alpha\beta}:= (\id_{F^*L}\otimes h_{\alpha\beta})\circ F^*u
\in
\operatorname{Hom}(F^*L^{-1},F^*L)|_{U_{\alpha\beta}}.
\]

We now use the exponential presentation of the inverse Cartier transform \cite[Section 2.2]{LSZ15}, with the sign convention specified below; compare \cite[Theorem 2.8]{OV07}. On each $U_\alpha$, equip $F^*E|_{U_\alpha}$ with the connection $\nabla_{h, \alpha}$ and use the transition 
functions $G^h_{\alpha\beta}$ defined by
\begin{equation}\label{eq:inverse-local-gluing}
	\begin{aligned}
		\nabla_{h,\alpha}
		&=
		\nabla^{\rm can}+\zeta_\alpha(F^*\phi),\\
		G^h_{\alpha\beta}
		&=
		\exp\bigl(h_{\alpha\beta}(F^*\phi)\bigr)
		\id+h_{\alpha\beta}(F^*\phi),
	\end{aligned}
\end{equation}
to glue the local data $(F^*E|_{U_\alpha}, \nabla_{h, \alpha})$ into a new flat bundle.
The exponential truncates after the linear term because $\phi^2=0$. With the overlap convention fixed above, the relation $dh_{\alpha\beta}=\zeta_\beta-\zeta_\alpha$ ensures that these local connections glue. We write
\[
(V,\nabla):=\mathrm{OV}_h(E,\phi)
\]
for the resulting flat bundle.

\begin{prop}
	The underlying vector bundle $V$ fits into an exact sequence
	\begin{equation}\label{exact:extension-h}
		0\longrightarrow F^*L \longrightarrow V \longrightarrow F^*L^{-1} \longrightarrow0.
	\end{equation}
	Its extension class is
	\begin{equation}\label{eq:main-extension-formula}
		[V]=F^*[E]+h\smile F^*u \quad\text{in }H^1(X,F^*L^2).
	\end{equation}
	The connections induced by $\nabla$ on $F^*L$ and $F^*L^{-1}$ are the canonical connections.
\end{prop}

\begin{proof}
	Choose local splittings of \eqref{eq:Higgs-extension}, and let $r_{\alpha\beta} \in\operatorname{Hom} (L^{-1},L)|_{U'_{\alpha\beta}}$ be the cocycle representing $[E]$. The induced local splittings of $F^*E$ have upper-right transition entries $F^*r_{\alpha\beta}$. Since $\theta$ annihilates $L$ and maps $L^{-1}$ into $L\otimes\Omega_{X'}$, the automorphism $G^h_{\alpha\beta}$ acts identically on the subbundle $F^*L$ and on the quotient $F^*L^{-1}$. Its upper-right entry is precisely $\nu^h_{\alpha\beta}$. Consequently, relative to the local splittings and the fixed end terms, the gluing maps for $V$ are
	\[
	\begin{pmatrix}
		1&F^*r_{\alpha\beta}+\nu^h_{\alpha\beta}\\
		0&1
	\end{pmatrix}.
	\]
	Therefore, by applying Proposition~\ref{extension-sum of classes} to the vector bundle $F^*E$ and the transition function $G^h_{\alpha\beta}$, we obtain \eqref{eq:main-extension-formula}.
	
	Finally, the correction term $\zeta_\alpha(F^*\phi)$ annihilates $F^*L$ and has image in $F^*L\otimes\Omega_X$. Thus its induced actions on the subbundle and quotient vanish, so both induced connections are canonical.
\end{proof}

\begin{rema}[Sign convention]\label{rem:p-curvature-sign}
	According to \cite{Ka70}, we use the convention $\psi_{\nabla}(D)	= \nabla_{h,D}^{\,p}-\nabla_{h,D^{[p]}}$. With this convention, the plus sign in \eqref{eq:inverse-local-gluing} gives $p$-curvature represented locally by $-F^*\phi$. In particular, the Cartier descent of the off-diagonal $p$-curvature is $-u$.
	
	This sign can be checked in an \'etale coordinate $z$. If
	\[
	\zeta(dz')=z^{p-1}dz, \qquad N^2=0,
	\]
	with $N$ horizontal for the canonical connection, then the connection $d+z^{p-1}N\,dz$ satisfies
	\[
	\psi\left(\frac{\partial}{\partial z}\right)=-N.
	\]
	Replacing $\phi$ by $-\phi$ gives the other common sign convention. Throughout this subsection, we retain the plus signs in \eqref{eq:inverse-local-gluing} and \eqref{eq:main-extension-formula}.
\end{rema}

The dependence on the choice of $h$ can also be expressed directly in terms of extension classes. According to Lemma \ref{lem:Fro-liftings}, let write $h'=h+F^*h_0$ with $h_0\in H^1(X',T_{X'})$. For the fixed Higgs bundle $(E,\phi)$, write
\[
(V_{h'},\nabla_{h'}) = \mathrm{OV}_{h'}(E,\phi).
\]
Then, equation (\ref{eq:main-extension-formula}) implies that
\[
[V_{h'}]= F^*[E]+h\smile F^*u+F^*(u\smile h_0).
\]
Let $E'$ be the extension
\[
0\longrightarrow L\xrightarrow{i_{E'}}E' \xrightarrow{\pi_{E'}}L^{-1}\longrightarrow0
\]
with class $[E']=[E]+u\smile h_0$ in $H^1(X',L^2)$, and equip it with the Higgs field 
\[
\phi' = (i_{E'}\otimes\id_{\Omega_{X'}})\circ u\circ \pi_{E'}.
\]
Then, by pulling back by $F$, we can write $[V_{h'}]$ as the underlying bundle of $\mathrm{OV}_h(E',\phi')$, and it is obvious to see that
\[
[V_{h'}]=F^*[E']+h\smile F^*u.
\]

\begin{prop}\label{prop:inverse-change-of-h}
	There is an isomorphism of flat bundles
	\[
	\mathrm{OV}_{h'}(E,\phi) \simeq \mathrm{OV}_h(E',\phi').
	\]
	Moreover, changing $h$ by $F^*h_0$ corresponds to adding $u\smile h_0$ to the extension class of $E$.
\end{prop}

\begin{proof}
	It will be enough to compare the connections on the vector bundle $V_{h'}$. We choose a 
    cocycle $h_{0,\alpha\beta}$ representing $h_0$ and represent $h'$ by
	\[
	h'_{\alpha\beta} = h_{\alpha\beta}+F^*h_{0,\alpha\beta}.
	\]
	Since $F^*h_{0,\alpha\beta}$ is horizontal, we have
	\[
	dh'_{\alpha\beta} = dh_{\alpha\beta} = \zeta_\beta-\zeta_\alpha.
	\]
	Hence the same local maps $\zeta_\alpha$ may be used for both $h$ and $h'$.

	If $[E]$ is represented by $r_{\alpha\beta}$, then $[E']$ is represented by
	\[
	r'_{\alpha\beta} = r_{\alpha\beta} + (\id_L\otimes h_{0,\alpha\beta})\circ u.
	\]
	It follows that the gluing maps for $\mathrm{OV}_h(E',\phi')$ coincide with those for $\mathrm{OV}_{h'}(E,\phi)$. Their local connection matrices also coincide, because $\phi$ and $\phi'$ have the same off-diagonal entry $u$. These local identifications therefore glue to the asserted isomorphism of flat bundles.
\end{proof}

\subsection{A normalized criterion for connections on extensions}

Fix a line bundle $L$ on $X'$, and consider an extension
\begin{equation}\label{eq:extension-for-criterion}
	0\longrightarrow F^*L\xrightarrow{i_V}V \xrightarrow{\pi_{V}}F^*L^{-1}\longrightarrow0.
\end{equation}
Its determinant is trivial, and its extension class is
\[
v=[V]\in \operatorname{Ext}^1_X(F^*L^{-1},F^*L) \simeq H^1(X,F^*L^2).
\]

Since $F$ is finite, $R^iF_*=0$ on quasi-coherent sheaves for $i>0$. Moreover, the projection formula, applied to $L^2$ and $\OO_X$, gives
\[
F_*F^*L^2 = F_*(F^*L^2\otimes\OO_X) \simeq L^2\otimes F_*\OO_X.
\]
Consequently,
\[
H^1(X,F^*L^2) \simeq H^1(X',F_*F^*L^2) \simeq H^1(X',L^2\otimes F_*\OO_X).
\]
Using these identifications, tensoring \eqref{ex:exact-functions} with $L^2$ induces the map
\[
d_{L^2}\colon H^1(X,F^*L^2) \longrightarrow H^1(X',L^2\otimes B).
\]

For a section $u\in H^0(X',L^2\otimes\Omega_{X'})$, the cup product followed by contraction defines
\[
[Ca]\smile u\in H^1(X',L^2\otimes B).
\]
Equivalently, this is the image of $u$ under the connecting homomorphism
\[
\delta_{L^2}\colon H^0(X',L^2\otimes\Omega_{X'}) \longrightarrow H^1(X',L^2\otimes B)
\]
associated with
\[
0\longrightarrow L^2\otimes B \longrightarrow L^2\otimes F_*\Omega_X \xrightarrow{\id_{L^2}\otimes C} L^2\otimes\Omega_{X'} \longrightarrow0.
\]
Indeed, viewing $u$ as a morphism $L^{-2}\to\Omega_{X'}$, the class $\delta_{L^2}(u)$ corresponds to the pullback of the Cartier extension along $u$, under the canonical identification 
\[
\operatorname{Ext}^1_{X'}(L^{-2},B) \simeq H^1(X',L^2\otimes B).
\]
This pullback is precisely the Yoneda product, so $\delta_{L^2}(u)=[Ca]\smile u$.

\begin{rema}[Cartier normalization] \label{rem:Cartier-normalization}
	A connection $\nabla$ on $V$ inducing the canonical connections on the subbundle and quotient of \eqref{eq:extension-for-criterion} has, in local splittings, the form
	\[
	\nabla_\alpha=
	\begin{pmatrix}
		\nabla^{\rm can}&\varphi_\alpha\\
		0&\nabla^{\rm can}
	\end{pmatrix}.
	\]
	By adjunction and the projection formula, we regard $\varphi_\alpha$ as a local section of $L^2\otimes F_*\Omega_X$ over $X'$. The compatibility of the local connections, together with $C\circ d=0$, implies that their images 
    by $\mathrm{id}_{L^2} \otimes C$ glue to a global section
	\[
	u_\nabla\in H^0(X',L^2\otimes\Omega_{X'}).
	\]
	This section is independent of the chosen local splittings. 
	
	For a prescribed section $u\in H^0(X',L^2\otimes\Omega_{X'})$, we say that the off-diagonal term is \emph{Cartier-normalized by $u$} if $u_\nabla=u$, or equivalently,
	\[
	(\id_{L^2}\otimes C)(\varphi_\alpha)=u
	\]
	on every member of the covering.
\end{rema}

\begin{prop}\label{prop:connection-criterion}
	Let $v\in H^1(X,F^*L^2)$ be the class of \eqref{eq:extension-for-criterion}, and fix $u\in H^0(X',L^2\otimes\Omega_{X'})$. Then $V$ admits a connection inducing the canonical connections on its subbundle and quotient, with off-diagonal term Cartier-normalized by $u$ in the sense of Remark \ref{rem:Cartier-normalization},
	if and only if 
	\begin{equation}\label{eq:connection-criterion}
		d_{L^2}(v)=[Ca]\smile u.
	\end{equation}
\end{prop}

\begin{proof}
	Choose a sufficiently fine affine covering $\{U'_\alpha\}$ of $X'$, and put $U_\alpha=F^{-1}(U'_\alpha)$. Choose local splittings
	\[
	s_\alpha\colon F^*L^{-1}|_{U_\alpha} \longrightarrow V|_{U_\alpha}
	\]
	of $\pi_V$, and write
	\[
	s_\beta-s_\alpha=i\circ v_{\alpha\beta}.
	\]
	Then $\{v_{\alpha\beta}\}$ represents $v$, where
	\[
		v_{\alpha\beta}	\in \operatorname{Hom}_{U_{\alpha\beta}}
		(F^*L^{-1},F^*L) \simeq \operatorname{Hom}_{U'_{\alpha\beta}}
		(L^{-2},F_*\OO_X).
	\]
	We use the same notation for these corresponding homomorphisms. Relative to the local splittings and the fixed end terms, the transition maps of $V$ are
	\[
	g_{\alpha\beta}	= 
	\begin{pmatrix}
		1&v_{\alpha\beta}\\
		0&1
	\end{pmatrix}.
	\]
	The class $d_{L^2}(v)$ is represented by
	\[
	dv_{\alpha\beta}:=d\circ v_{\alpha\beta}\in \operatorname{Hom}_{U'_{\alpha\beta}}(L^{-2},B)	= H^0(U'_{\alpha\beta},L^2\otimes B).
	\]
	Under adjunction, this is the covariant derivative of $v_{\alpha\beta}$ for the canonical connections on the two end terms.

    The pullback of the Cartier extension along \(u\colon L^{-2}\to\Omega_{X'}\) represents \([Ca]\smile u\) and fits into the following commutative diagram:
	\[
	\begin{tikzcd}
		0 \arrow[r] &
		B \arrow[r] \arrow[d,equal] &
		\mathcal E_{\mathrm{Ca},u} \arrow[r,"q_u"] \arrow[d,"\rho_u"] &
		L^{-2} \arrow[r] \arrow[d,"u"] &
		0\\
		0 \arrow[r] &
		B \arrow[r] &
		F_*\Omega_X \arrow[r,"C"] &
		\Omega_{X'} \arrow[r] &
		0.
	\end{tikzcd}
	\]

	Choose local Cartier splittings
	\[
	\zeta_\alpha\colon
	\Omega_{X'}|_{U'_\alpha}	\longrightarrow F_*\Omega_X|_{U'_\alpha}, \qquad C\circ\zeta_\alpha=\id,
	\]
	and put $\zeta_\alpha^u:=\zeta_\alpha\circ u$. Then $C\circ\zeta_\alpha^u=u$, and the maps
	\[
	a\longmapsto\bigl(\zeta_\alpha^u(a),a\bigr)
	\]
	are local splittings of $q_u$. It follows that
	\[
	[Ca]\smile u = \bigl[\{\zeta_\beta^u-\zeta_\alpha^u\}_{\alpha\beta}\bigr].
	\]
	
	Suppose first that
	\[
	d_{L^2}(v)=[Ca]\smile u.
	\]
	This equality of cohomology classes means that there is a $0$-cochain
	\[
	b_\alpha\in \operatorname{Hom}_{U'_\alpha}(L^{-2},B)
	\]
	such that
	\[
	dv_{\alpha\beta} - (\zeta_\beta^u-\zeta_\alpha^u) = b_\beta-b_\alpha.
	\]
	Define $\varphi_\alpha:=\zeta_\alpha^u+b_\alpha$. Since $b_\alpha$ takes values in $B=\ker C$, we obtain
	\[
	C\circ\varphi_\alpha=u, \qquad dv_{\alpha\beta} = \varphi_\beta-\varphi_\alpha.
	\]
	By adjunction, these maps correspond to
	\[
	\begin{aligned}
		\varphi_\alpha
		&\in
		\operatorname{Hom}_{U'_\alpha}
		(L^{-2},F_*\Omega_X) \simeq
		\operatorname{Hom}_{U_\alpha}
		(F^*L^{-1},F^*L\otimes\Omega_X).
	\end{aligned}
	\]
	We therefore equip $(F^*L\oplus F^*L^{-1})|_{U_\alpha}$ with
	\[
	\nabla_\alpha = 
	\begin{pmatrix}
		\nabla^{\rm can}&\varphi_\alpha\\
		0&\nabla^{\rm can}
	\end{pmatrix}.
	\]
	
	We verify that these local connections glue to a global connection. Let $\nabla_0$ be the diagonal canonical connection on $F^*L\oplus F^*L^{-1}$, and write
	\[
	\nabla_\alpha=\nabla_0+A_\alpha, \qquad 
	A_\alpha=
	\begin{pmatrix}
		0&\varphi_\alpha\\
		0&0
	\end{pmatrix}.
	\]
	Writing $d$ for the differential induced by $\nabla_0$ on endomorphisms, the connection transformation law is 
	\[
	A_\beta	= g_{\alpha\beta}^{-1}A_\alpha g_{\alpha\beta} + g_{\alpha\beta}^{-1}dg_{\alpha\beta}.
	\]
	Since
	\[
	g_{\alpha\beta}^{-1} =
	\begin{pmatrix}
		1&-v_{\alpha\beta}\\
		0&1
	\end{pmatrix},
	\]
	direct multiplication gives
	\[
	g_{\alpha\beta}^{-1}A_\alpha g_{\alpha\beta} + g_{\alpha\beta}^{-1}dg_{\alpha\beta}	=
	\begin{pmatrix}
		0&\varphi_\alpha+dv_{\alpha\beta}\\
		0&0
	\end{pmatrix}.
	\]
	Thus the local connections glue precisely when
	\begin{equation}\label{eq:connection-gluing}
		dv_{\alpha\beta} = \varphi_\beta-\varphi_\alpha.
	\end{equation}
	This relation holds by construction. The resulting global connection on $V$ induces the canonical connections on its subbundle and quotient. Moreover, $C\circ\varphi_\alpha=u$, so its off-diagonal term is Cartier-normalized by $u$. Conversely, suppose that $V$ admits a connection with the stated properties.
	In the chosen local splittings, write
	\[
	\nabla_\alpha = 
	\begin{pmatrix}
		\nabla^{\rm can}&\varphi_\alpha\\
		0&\nabla^{\rm can}
	\end{pmatrix},
	\]
	viewing $\varphi_\alpha$ by adjunction as a map $L^{-2}\to F_*\Omega_X$. The compatibility of these local connections gives
	\[
	dv_{\alpha\beta} = \varphi_\beta-\varphi_\alpha,
	\]
	while the prescribed Cartier normalization gives
	\[
	C\circ\varphi_\alpha=u.
	\]
	Since also $C\circ\zeta_\alpha^u=u$, the differences
	\[
	b_\alpha:=\varphi_\alpha-\zeta_\alpha^u
	\]
	take values in $B$. Hence,
	\[
	dv_{\alpha\beta} = (\zeta_\beta^u-\zeta_\alpha^u)+(b_\beta-b_\alpha).
	\]
	The cocycles $\{dv_{\alpha\beta}\}$ and $\{\zeta_\beta^u-\zeta_\alpha^u\}$ therefore represent the same class of cohomology, proving
	\[
	d_{L^2}(v)=[Ca]\smile u.
	\]
	
	Finally, the constructed connection is integrable because $X$ is a curve. Its induced determinant connection is the trivial connection on $F^*L\otimes F^*L^{-1}\simeq\OO_X$.
\end{proof}

\begin{rema}\label{rem:sign-p-curvature}
    Under the $p$-curvature convention of Remark~\ref{rem:p-curvature-sign}, a connection Cartier-normalized by $u$ has off-diagonal $p$-curvature $-F^*u$.
\end{rema}

When $E$ has arbitrary determinant, a similar argument yields the following result.
\begin{rema}\label{rema:not-trivialdeteminant-connection-criterion}
	Let $0\longrightarrow F^*L\xrightarrow{i_V}V \xrightarrow{\pi_V}F^*Q\longrightarrow0$ be an extension of $F^*Q$ by $F^*L$ with extension class $v\in H^1(X,F^*(LQ^{-1}))$, and fix $u\in H^0(X',LQ^{-1}\otimes\Omega_{X'})$. Then $V$ admits a connection inducing the canonical connections on $F^*L$ and $F^*Q$, with off-diagonal term Cartier-normalized by $u$ in the sense of Remark~\ref{rem:Cartier-normalization}, if and only if 
\begin{equation*}%\label{eq:connection-criterion}
	d_{L^2}(v)=[Ca]\smile u,
\end{equation*}
where $d\colon H^1(F_*\OO\otimes Q^{-1}L) \longrightarrow H^1(B\otimes Q^{-1}L)$ is the map induced by $d\colon F_*\OO\to B$.

Indeed, the proof of Proposition~\ref{prop:connection-criterion} applies with $LQ^{-1}$ in place of $L^2$. Thus the same criterion applies to extensions associated with a Higgs bundle
\[
0\longrightarrow L\longrightarrow E\longrightarrow Q \longrightarrow0
\]
without requiring $\det E=L\otimes Q$ to be trivial.
\end{rema}

\subsection{The Cartier transform}

Fix $h\in\LL(X)$ with local data $(h_{\alpha\beta},\zeta_\alpha)$ as above. We describe the Cartier transform $C_{\exp}$ in terms of extension classes, following \cite[Section~2.3]{LSZ15}.

Let $(V,\nabla)$ be a rank-$2$ flat bundle on $X$ satisfying
\[
\det V\simeq\OO_X, \qquad \psi_\nabla\ne0, \qquad \psi_\nabla^2=0.
\]
Denote by $K=\ker\psi_\nabla$. Then $K$ is a line subbundle and $V/K\simeq K^{-1}$. The $p$-curvature factors through the quotient and takes values in $K\otimes F^*\Omega_{X'}$, giving the commutative diagram with exact rows
\[
\begin{tikzcd}[column sep=small,row sep=large]
	0\arrow[r]& 
	K\arrow[r]&
	V\arrow[r]\arrow[d,"\psi_\nabla"]&
	K^{-1}\arrow[r]
	\arrow[dll,bend right=45,"\beta"']&
	0\\
	0\arrow[r]&
	K\otimes F^*\Omega_{X'}\arrow[r]&
	V\otimes F^*\Omega_{X'}\arrow[r]&
	K^{-1}\otimes F^*\Omega_{X'}\arrow[r]&
	0.
\end{tikzcd}
\]

By \cite[Proposition~5.2, (5.2.3)]{Ka70}, the $p$-curvature is horizontal as a map
\[
\psi_\nabla\colon(V,\nabla)\longrightarrow \bigl(V\otimes F^*\Omega_{X'}, \nabla\otimes\nabla^{\rm can}\bigr).
\]
Consequently, $\nabla$ preserves $K$ and induces a connection $\overline\nabla$ on $V/K$. Both induced connections have vanishing $p$-curvature. Cartier descent therefore gives line bundles $L,Q$ on $X'$ with canonical connections
\[
(K,\nabla|_K)\simeq(F^*L,\nabla^{\rm can}), \qquad (V/K,\overline\nabla)\simeq(F^*Q,\nabla^{\rm can}).
\]
Thus we obtain an extension
\begin{equation}\label{eq:flat-extension-recovery}
	0\longrightarrow F^*L\longrightarrow V \longrightarrow F^*Q\longrightarrow0.
\end{equation}
Moreover, the determinant assumption gives $F^*(L\otimes Q)\simeq\OO_X$.

The induced map
\[
\beta\colon F^*Q\longrightarrow F^*L\otimes F^*\Omega_{X'}
\]
is horizontal and hence descends under Cartier descent. With the sign convention of Remark~\ref{rem:p-curvature-sign}, there exists $0\ne u\in H^0(X',LQ^{-1}\otimes\Omega_{X'})$ such that
\[
\beta=-F^*u.
\]
Equivalently, $u$ is the Cartier-normalized off-diagonal section associated with $\nabla$.

Write
\[
v=[V]\in \Ext^1_X(F^*Q,F^*L) \simeq H^1(X,F^*(LQ^{-1})) \simeq H^1(X',LQ^{-1}\otimes F_*\OO_X).
\]
We also use $d$ to denote the cohomology maps induced by $d: F_*\OO_X\to B$; in particular,
\[
d\colon H^1(X',LQ^{-1}\otimes F_*\OO_X) \longrightarrow H^1(X',LQ^{-1}\otimes B).
\]

\begin{prop}[Recovery of the Higgs extension]\label{prop:recover-Higgs-extension}
	There exists an extension class $e\in H^1(X',LQ^{-1})$ defining an extension $0\longrightarrow L\xrightarrow{i_E}E \xrightarrow{\pi_E}Q\longrightarrow0$ such that
	\begin{equation}\label{eq:recovered-extension-class}
		v=F^*e+(h\smile F^*u).
	\end{equation}
	Moreover, $e$ can be chosen so that the Higgs field $\phi=	(i_E\otimes\id_{\Omega_{X'}})\circ u\circ \pi_E$ satisfies
	\[
	\mathrm{OV}_h(E,\phi)\simeq(V,\nabla).
	\]
\end{prop}

\begin{proof}
	Applying Proposition \ref{prop:connection-criterion} in the generality of Remark~\ref{rema:not-trivialdeteminant-connection-criterion}, we obtain
	\[
	d(v)=[Ca]\smile u \in H^1(X',LQ^{-1}\otimes B).
	\]
	
	The section $u\in H^0(X',LQ^{-1}\otimes\Omega_{X'})$, viewed as a morphism $T_{X'}\to LQ^{-1}$, also induces a map by taking the cup product
	\[
	{-}\smile u\colon
	H^1(X',T_{X'}\otimes B)	\longrightarrow H^1(X',LQ^{-1}\otimes B).
	\]
	Hence, $[Ca]\smile u$  belongs to $H^1(X',LQ^{-1}\otimes B)$, where $[Ca]$ is the Cartier class.

	Under the canonical identification
	\[
	H^1(X,F^*T_{X'})
	\simeq H^1(X',T_{X'}\otimes F_*\OO_X),
	\]
	we use the same notation $h$ for its image. By definition, $d(h)=[Ca]$. W clearly have the commutative
	diagram
	\[
	\begin{tikzcd}[column sep=large,row sep=large]
		H^1(X',T_{X'}\otimes F_*\OO_X) 	\arrow[r,"d"] \arrow[d,"{-}\smile u"']& H^1(X',T_{X'}\otimes B) \arrow[d,"{-}\smile u"]\\
		H^1(X',LQ^{-1}\otimes F_*\OO_X) \arrow[r,"d"]& H^1(X',LQ^{-1}\otimes B),
	\end{tikzcd}
	\]
	where all cohomology groups are taken over $X'$. It follows that
	\[
	d([h]\smile u)=[Ca]\smile u.
	\]
	Hence, we have two extension classes $h\smile u$ and $v$ lying over $[Ca]\smile u$ under the map $d: H^1(LQ^{-1}\otimes F_*\OO_X)\longrightarrow H^1(LQ^{-1}\otimes B)$. This shows that
	\[
	v-(h\smile u)\in\ker d.
	\]

	Tensoring \eqref{ex:exact-functions} with $LQ^{-1}$ and taking cohomology gives
	\[
	\ker d = \operatorname{im}\!\left(F^*\colon H^1(X',LQ^{-1}) \longrightarrow H^1(X,F^*(LQ^{-1})) \right).
	\]
	Consequently, there exists a class $e$ with
	\[
	F^*e=v-(h\smile F^*u).
	\]
	
	To obtain an extension compatible with the given connection, apply the exponential twisting construction of \cite[Section~2.3]{LSZ15}. Locally, this replaces the connection and gluing data by
	\[
	\nabla'_\alpha = \nabla_\alpha+\zeta_\alpha(\psi_\nabla),
	\qquad
	J_{\alpha\beta} = \exp\!\bigl(h_{\alpha\beta}(\psi_\nabla)\bigr).
	\]
	These data define a flat bundle $(H',\nabla')$ with vanishing $p$-curvature. Since $\beta=-F^*u$, its underlying extension class is
	\[
	[H']=v-(h\smile F^*u).
	\]
	The construction leaves the induced connections on the subbundle and quotient unchanged. Cartier descent therefore produces an extension
	\[
	0\longrightarrow L\longrightarrow E \longrightarrow Q\longrightarrow0
	\]
	whose class $e$ satisfies \eqref{eq:recovered-extension-class}. The map $u$ defines the Higgs field $\phi$. The inverse exponential construction recovers $(V,\nabla)$. Hence, $$\mathrm{OV}_h(E,\phi)\simeq(V,\nabla).$$
\end{proof}

\begin{rema}[The determinant under Cartier descent]\label{rem:determinant-Cartier-descent}
	The exponential twisting construction preserves the determinant connection: its connection corrections have trace zero, and its gluing corrections have determinant one. Consequently,
	\[
	\bigl(F^*(\det E),\nabla^{\rm can}\bigr) \simeq (\det V,\det\nabla).
	\]
	In particular,
	\[
	\det E=L\otimes Q, \qquad F^*(\det E)\simeq\OO_X.
	\]
	This does not imply $\det E\simeq\OO_{X'}$, since Frobenius pullback on line bundles need not be injective. Moreover,
	\[
	\psi_{\det\nabla} = \operatorname{tr}(\psi_\nabla)=0
	\]
	already follows from the nilpotence of $\psi_\nabla$; this vanishing alone does not force $\det E$ to be trivial.
\end{rema}

\begin{coro}\label{cor:recovery-trivial-determinant}
	Under the preceding assumptions, suppose in addition that
	\[
	(\det V,\det\nabla)\simeq(\OO_X,d).
	\]
	Then $Q\simeq L^{-1}$, and the construction in Proposition~\ref{prop:recover-Higgs-extension} gives an extension
	\[
	0\longrightarrow L\longrightarrow E \longrightarrow L^{-1}\longrightarrow0
	\]
	with $\det E\simeq\OO_{X'}$. Its class $e\in H^1(X',L^2)$ satisfies
	\[
	v=F^*e+h\smile F^*u, \qquad \mathrm{OV}_h(E,\phi)\simeq(V,\nabla).
	\]
\end{coro}

\begin{proof}
	The horizontal determinant trivialization identifies the quotient flat line bundle with the dual of $K$:
	\[
	(V/K,\overline\nabla) \simeq \bigl(K^\vee,(\nabla|_K)^\vee\bigr).
	\]
	Cartier descent therefore gives $Q\simeq L^{-1}$. The assertions now follow from Proposition~\ref{prop:recover-Higgs-extension}, with $\det E\simeq L\otimes L^{-1}\simeq\OO_{X'}$.
\end{proof}

% Sources 9--12 supply the applications in Section 3.
\section{Stability of the underlying bundles} \label{sec:semistability}

\subsection{Semistability of pairs and of their underlying bundles}

Fix $h\in\LL(X)$, and let $(E,\phi)$ be a rank-$2$ nilpotent Higgs bundle on $X'$ with $\det E\simeq\OO_{X'}$. Write
\[
(V,\nabla)=\mathrm{OV}_h(E,\phi).
\]
The preceding constructions give $\det V\simeq\OO_X$, with the induced trivial determinant connection. In particular,
\[
\mu(E)=\mu(V)=0,
\]
where $\mu(W)=\deg(W)/\operatorname{rk}(W)$. Throughout this subsection, a line subbundle means a saturated rank-one subsheaf. In this setting, $(E,\phi)$ is semistable if $\deg A\le0$ for every $\phi$-invariant line subbundle $A\subset E$. Similarly, $(V,\nabla)$ is semistable if $\deg B\le0$ for every $\nabla$-invariant line subbundle $B\subset V$. Strict inequalities define stability. For the underlying vector bundles, the corresponding inequalities must hold for all line subbundles. Thus semistability of an underlying bundle implies semistability of the pair.

Suppose first that $\phi \ne0$, and let $L=\ker\phi$ as in \eqref{eq:Higgs-extension}. Then $L$ is the unique $\phi$-invariant line subbundle of $E$. Indeed, after trivializing $\Omega_{X'}$ at the generic point, $\phi$ becomes a nonzero square-zero endomorphism of a two-dimensional vector space. Its only invariant line is its kernel. An invariant line subbundle therefore agrees generically with $L$, and saturation implies that it equals $L$. Consequently,
\begin{equation}\label{eq:Higgs-semistability}
	\begin{aligned}
		(E,\phi)\text{ is semistable} 	&\quad\Longleftrightarrow\quad \deg L\le0,\\
		(E,\phi)\text{ is stable} 		&\quad\Longleftrightarrow\quad \deg L<0.
	\end{aligned}
\end{equation}

The preservation statement below follows from \cite[Theorem 4.7 and Remark 4.10]{dCGZ24}. 
\begin{prop}\label{prop:pair-semistability}
	The correspondence $\mathrm{OV}_h$ preserves semistability and stability of rank-$2$ nilpotent Higgs bundles and flat bundles with trivial determinant.
	More precisely,
	\[
	\begin{aligned}
		(E,\phi)\text{ is semistable} &\quad\Longleftrightarrow\quad (V,\nabla)\text{ is semistable},\\
		(E,\phi)\text{ is stable} 	&\quad\Longleftrightarrow\quad (V,\nabla)\text{ is stable}.
	\end{aligned}
	\]\qed
\end{prop}

\begin{rema}
	Proposition~\ref{prop:pair-semistability} concerns 	Higgs and flat pairs. For example, if $E$ is semistable, then $\mathrm{OV}_h(E,0)=(F^*E,\nabla^{\rm can})$ is semistable as a flat bundle, although Frobenius pullback can destabilize the underlying vector bundle $F^*E$.
\end{rema}

\subsection{Stability of the underlying bundles}\label{subsec:H-stability}%not finish in the switing
Recall that the affine space $\mathcal L(X)$ is modeled on the vector space $H^1(X', T_{X'})$, 
so $\dim\mathcal L(X)=3g-3$. We recall from Definition \ref{def:lift-space} that 
$\widehat{\mathcal{L}(X)} := \pi^{-1}(\Bbbk[Ca]) \subset H^1(X,F^*T_{X'})$. Then these spaces satisfy $\PP(\widehat{\mathcal{L}(X)})=\LL(X)\sqcup\PP H^1(X',T_{X'})$, and hence
\begin{equation}\label{eq:projective-compactification}
 \PP(\widehat{\mathcal{L}(X)})\simeq\PP^{3g-3}.
\end{equation}

Consider the forgetful map
\[
\begin{tikzcd}
	\mathcal L(X)\subset\PP^{3g-3}=\PP(\widehat{\mathcal{L}(X)}) \arrow[r,dashed,"\overline\Phi"]&M_X(2,F^*\Omega_{X'})\\
	\PP^{3g-4}=\PP H^1(X',T_{X'})\arrow[u,hook]\arrow[ur,dashed]&
\end{tikzcd}
\]
given by
\[
\overline\Phi\bigl(0\longrightarrow\OO_X\longrightarrow \HH_h^2 \longrightarrow F^*\Omega_{X'}\longrightarrow0\bigr)
= \HH_h^2.
\]
and $\HH_h^2$ is the extension \eqref{eq:H2} defined by $h$.
We want to show that for a general class $h\in \mathcal L(X)$ the associated rank-$2$ vector bundle $\HH_h^2$ is stable. Since the stability condition is open, it will be enough to show that there exists an $e=[E]\in H^1(X', T _ {X'}) $ such that $F^*E\in H^1(X', F^*T _ {X '}) $ is stable.

Assume that $p \ne 2$. We consider a connected étale double cover $\pi: Y \to X$ and denote by 
$\sigma : Y \to Y$ the covering involution and by $\alpha \in \Pic^0(X)[2]$ its associated $2$-torsion 
line bundle, i.e. we have $\pi_*\mathcal{O}_Y = \mathcal{O}_X \oplus \alpha$. For every $M\in\Pic^0(Y)$, the bundle
$\pi_*M$ has rank $2$ and degree $0$. According to \cite[Corollaire 1.5.2 (1.5.2.1)]{EGAII}, we have
\[
\pi^*\pi_*M\simeq M\oplus\sigma^*M.
\]
Consequently, $\pi_*M$ is semistable.

Suppose that $\pi_*M$ is strictly semistable. It then contains a line subbundle $N\in\Pic^0(X)$. By adjunction, the inclusion $N\hookrightarrow\pi_*M$ induces a nonzero morphism $\pi^*N\longrightarrow M$ which is an isomorphism. Conversely, if $M\simeq\pi^*N$ for some $N\in\Pic^0(X)$, the projection formula gives
\[
\pi_*M\simeq N\otimes\pi_*\OO_{\widetilde X},
\]
which contains $N$ as a degree-$0$ line subbundle. Thus
\[
\pi_*M\text{ is strictly semistable} \quad\Longleftrightarrow\quad M\simeq\pi^*N \quad\text{for some }N\in\Pic^0(X).
\]

Moreover, $M$ descends to $X$ if and only if $\sigma^*M\simeq M$. Indeed, over the algebraically closed field $\Bbbk$, an isomorphism $\sigma^*M\simeq M$ can be rescaled to satisfy the cocycle condition, and hence defines a descent datum. We therefore obtain
\[
\pi_*M\text{ is strictly semistable} \quad\Longleftrightarrow\quad \sigma^*M\simeq M.
\]

Let $M\in \operatorname{Prym}(Y/X)$. The norm identity gives $M\otimes\sigma^*M \simeq\pi^*\operatorname{Nm}_{\pi}(M) \simeq\OO_{Y}$, so that $\sigma^*M\simeq M^{-1}$. 
Consequently, $\pi_*M$ is strictly semistable if and only if $M^{\otimes2}\simeq\OO_{Y}$, and hence it is equivalent to $M\in \operatorname{Prym}(Y/X)[2]$, where $\operatorname{Prym}(Y/X)[2]$ is the finite étale subgroup scheme of $2$-torsion points.

Since $\pi$ is étale, we have a Cartesian diagram
\begin{equation} \label{cdFrob}
\begin{tikzcd} 
 Y \ar[r,"F_Y"] \ar[d, "\pi"] &  Y' \ar[d, "\pi'"] \\
 X \ar[r, "F"] & X',
 \end{tikzcd}
\end{equation}
where the horizontal maps are the relative Frobenius maps. Let $\alpha' \in \mathrm{Pic}(X')[2]$ be the $2$-torsion
line bundle associated with $\pi'$.

\begin{lemm}
	Let $L\in\operatorname{Nm}_{\pi'}^{-1}(\alpha')$ and let $E_L:=\pi'_*L$. Then, both $F^*E_L$ and $E_L$ are stable whenever
	\[
	L\in \operatorname{Nm}_{\pi'}^{-1}(\alpha') 	\setminus \Pic^0(Y')[2p].
	\]
\end{lemm}
\begin{proof}
The determinant formula for the norm gives
\[
\det E_L \simeq \operatorname{Nm}_{\pi'}(L) \otimes\det(\pi'_*\OO_{Y'}) \simeq\alpha'^{\otimes2} \simeq\OO_{X'}.
\]
The Frobenius square \eqref{cdFrob} associated with $\pi$ and $\pi'$ is Cartesian, and hence $F^*E_L\simeq\pi_*F_{Y}^*L$. Compatibility of the norm with base change gives 
$\operatorname{Nm}_{\pi}(F_{Y}^*L)\simeq F^*\alpha'= \alpha^p = \alpha$, while
\[
\pi^*F^*\alpha' \simeq F_{Y}^*\pi'^*\alpha' \simeq\OO_{Y}.
\]
It follows that
\[
F_{Y}^*L\otimes\sigma^*F_{Y}^*L \simeq\pi^*\operatorname{Nm}_{\pi}(F_{Y}^*L) \simeq\OO_{Y}.
\]
Thus $\sigma^*F_{Y}^*L\simeq F_{Y}^*L^{-1}$, and the preceding
criterion yields
\[
\begin{aligned}
	F^*E_L\text{ is strictly semistable}&\quad
    \Longleftrightarrow\quad
	\sigma^*F_{Y}^*L\simeq F_{Y}^*L \quad
    \Longleftrightarrow\quad
	F_{Y}^*L^{\otimes2}\simeq\OO_{Y}.
\end{aligned}
\]
By the identification of Frobenius pullback on the Jacobian with Verschiebung \cite[Proposition~2.1, p.~79]{Oda69} and the Frobenius--Verschiebung relation, we have
\[
F_{Y}^*(L^{\otimes2})\simeq\OO_{Y} \quad\Longleftrightarrow\quad L^{\otimes2p}\simeq\OO_{Y'}.
\]

Therefore,
	\[
	F^*E_L\text{ is stable} \quad    \Longleftrightarrow\quad  L\notin\Pic^0(Y')[2p].
	\]
In particular, both $F^*E_L$ and $E_L$ are stable whenever
	\[
	L\in \operatorname{Nm}_{\pi'}^{-1}(\alpha') \setminus \Pic^0(Y')[2p].
	\]
\end{proof}

\begin{prop}
Assume that $g \geq 3$. A general extension class $e \in H^1(X', T_{X'}) = \mathrm{Ext}^1(\Omega_{X'}, \OO_{X'})$   
defines a stable bundle $E_e$ fitting in the exact sequence
\[
0\longrightarrow\OO_{X'}\longrightarrow E_e\longrightarrow\Omega_{X'} \longrightarrow 0,
\]
and such that $F^*E_e$ is stable.
\end{prop}

\begin{proof}
Fix a theta characteristic $\theta$ on $X'$ and consider the divisor in 
$\operatorname{Nm}_{\pi'}^{-1}(\alpha')$ defined by $\theta$:
\[
 \mathcal D_\theta=\{L\in \operatorname{Nm}_{\pi'}^{-1}(\alpha')\mid H^0(X',E_L\otimes\theta)\ne0\}.
\]
Since $\operatorname{Pic}^0(Y')[2p]$ is finite and since we assume $g \geq 3$, we can choose
\[
L\in \mathcal D_\theta\setminus \operatorname{Nm}_{\pi'}^{-1}(\alpha')[2p]
\]
Then $\pi_*(F_{Y}^*L)=F^*(\pi'_*L) = F^* E_L$ is stable. Thus, the rational map induced by pull-back by Frobenius 
$E \mapsto F^*E$ is defined on an open subset of the generalized Theta divisor 
$$ \Theta \subset \mathcal{M}(X',2,\Omega_{X'}),$$
since we have constructed a stable bundle $E = E_L \otimes \theta \in \Theta$ with $F^* E$ stable. It is
well-known that the rational classifying map $\mathbb{P} H^1(X',T_{X'}) \dashrightarrow \Theta$ is dominant, from
which we can conclude the proposition.
\end{proof}

Since stability is an open condition, we obtain

\begin{coro}
For a general $h \in \widehat{\mathcal{L}(X)}$ the rank-$2$ vector bundle $\HH_h^2$ is stable.
\end{coro}

\section{Very stability and finiteness of the restricted de Rham--Hitchin morphism}

\subsection{Homogeneity of the restricted Hodge--Hitchin morphism for \texorpdfstring{$t$}{t}-connections}

Assume $p>2$. Let $V\in\mathcal M(X,2)$ be a stable rank-$2$ vector bundle, and fix a trivialization $\det V\simeq\OO_X$.

\begin{defi}
	For $t\in\Bbbk$, a $t$-connection on $V$ is a $\Bbbk$-linear morphism $\nabla_t\colon V\longrightarrow\Omega_X^1\otimes V$ satisfying
	\[
	\nabla_t(fv)=f\nabla_t(v)+t\,df\otimes v
	\]
	for local sections $f$ of $\OO_X$ and $v$ of $V$.
\end{defi}

A $0$-connection is a Higgs field, whereas a $1$-connection is an ordinary connection. We recall how the sheaf of the principal parts gives a vector-space description of connections with varying $t$ (cf.~\cite[Lemmas~3.1--3.2]{LP01}).

Let
\[
\mathcal P_X^1 = p_{1*}\bigl(\OO_{X\times X}/\mathcal I_\Delta^2\bigr),
\]
where $p_1$ is the first projection and $\mathcal I_\Delta$ is the ideal of the diagonal. The two projections give $\mathcal P_X^1$ its left and right $\OO_X$-module structures. There is a canonical exact sequence
\begin{equation}\label{exact:principal-part}
	0\longrightarrow\Omega_X^1 \longrightarrow\mathcal P_X^1 \xrightarrow{\rho}\OO_X \longrightarrow0,
\end{equation}
in which $df$ is represented locally by $1\otimes f-f\otimes1$. Tensoring with $V$ using the right module structure gives an exact sequence of left $\OO_X$-modules
\[
\begin{tikzcd}
	0\arrow[r]& \Omega_X^1\otimes V\arrow[r]& \mathcal P_X^1(V)\arrow[r,"\rho_V"]& V\arrow[r]&0,
\end{tikzcd}
\]
where $\mathcal P_X^1(V):=\mathcal P_X^1\otimes_{\OO_X}V$.

Composition with $\rho_V$ defines
\[
\pi_0\colon \Hom_{\OO_X}\bigl(V,\mathcal P_X^1(V)\bigr) \longrightarrow \Hom_{\OO_X}(V,V).
\]
Set
\[
\mathrm{tConn}(V) := \pi_0^{-1}(\Bbbk\cdot\id_V).
\]
This is a vector space, whose interpretation in terms of $t$-connections is given by the following lemma.

\begin{lemm}\label{lem:t-connections-principal-parts}
	For each $t\in\Bbbk$, the formula
	\[
	\nabla_t(v)=(t\otimes1)\otimes v-\sigma(v)
	\]
	gives a bijection between $\OO_X$-linear maps $\sigma\colon V\to\mathcal P_X^1(V)$ satisfying $\rho_V\circ\sigma=t\,\id_V$ and $t$-connections on $V$.
\end{lemm}

\begin{proof}
	Suppose first that $\rho_V\circ\sigma=t\,\id_V$. Then
	\[
	\rho_V\bigl((t\otimes1)\otimes v-\sigma(v)\bigr)=tv-tv=0,
	\]
	so $\nabla_t(v)$ lies in $\Omega_X^1\otimes V$. Using the right $\mathcal{O}_X$-module structure in the tensor product and the identity $df=1\otimes f-f\otimes1$, we obtain
	\[
	\begin{aligned}
		\nabla_t(fv)&=(t\otimes f)\otimes v-f\sigma(v)\\
		&=f\bigl((t\otimes1)\otimes v-\sigma(v)\bigr)+t\,df\otimes v\\
		&=f\nabla_t(v)+t\,df\otimes v.
	\end{aligned}
	\]
	Thus $\nabla_t$ is a $t$-connection.
	
	Conversely, given a $t$-connection $\nabla_t$, define
	\[
	\sigma(v)=(t\otimes1)\otimes v-\nabla_t(v),
	\]
	viewing the differential-form term inside $\mathcal P_X^1(V)$. The same identity and the Leibniz rule give
	\[
	\begin{aligned}
		\sigma(fv)
		&=f\bigl((t\otimes1)\otimes v\bigr)
		+t\,df\otimes v
		-f\nabla_t(v)-t\,df\otimes v=f\sigma(v).
	\end{aligned}
	\]
	Hence $\sigma$ is $\OO_X$-linear. Moreover, since $\rho_V$ vanishes on $\Omega_X^1\otimes V$, we obtain
	\[
	(\rho_V\circ\sigma)(v) = \rho_V\bigl((t\otimes1)\otimes v\bigr) - \rho_V\bigl(\nabla_t(v)\bigr) = tv.
	\]
	The two constructions are inverse to one another.
\end{proof}

We henceforth identify $\mathrm{tConn}(V)$ with the vector space of $t$-connections on $V$, allowing $t$ to vary. Since $V$ is stable of degree zero, it admits an ordinary connection \cite{Bi05}. Choose one, denoted by $\nabla_0$. If $\det\nabla_0=d+\eta$, replacing $\nabla_0$ by $\nabla_0-\frac12\eta\,\id_V$ ensures that
\[
\det\nabla_0=d.
\]
The sequence
\[
0\longrightarrow H^0(X,\End(V)\otimes\Omega_X^1) \longrightarrow \mathrm{tConn}(V) \xrightarrow{\nabla_t\mapsto t} \Bbbk \longrightarrow0
\]
is then split by $t\mapsto t\nabla_0$. Set
\[
\mathrm{tConn}^0(V) := \{\nabla_t\in\mathrm{tConn}(V)\mid\det\nabla_t=t\,d\}.
\]
Writing $\End^0(V)$ for trace-free endomorphisms, the unique decomposition $\nabla_t=t\nabla_0+\Phi$ identifies
\[
\mathrm{tConn}^0(V)\simeq H^0(X,\End^0(V)\otimes\Omega_X^1)\oplus \Bbbk.
\]
Because of stability, and using Riemann--Roch and Serre duality, we obtain
\begin{equation}\label{eq:dim-tConn}
	n:=h^0(X,\End_0(V)\otimes\Omega_X^1)=3g-3.
\end{equation}
Choosing a basis $\Phi_1,\ldots,\Phi_n$ therefore gives coordinates $x_1,\ldots,x_n,t$ through
\[
\nabla_t=t\nabla_0+\sum_{j=1}^n x_j\Phi_j.
\]

For a local vector field $D$, the $p$-curvature of $\nabla_t$ is
\[
\psi_{\nabla_t}(D) = \nabla_t(D)^p-t^{p-1}\nabla_t(D^{[p]}),
\]
where $D^{[p]}=D^p$ is the restricted $p$-th power. It determines a section of $\End(V)\otimes F^*\Omega_{X'}^1$. For $\nabla_t\in\mathrm{tConn}^0(V)$, we have  $\tr\bigl(\psi_{\nabla_t}\bigr) = \psi(\det\nabla_t) = 0$. Thus its characteristic polynomial is determined by its determinant, which descends under Frobenius to a unique section
\[
a_2(\nabla_t)\in \mathcal A_{X'}^0:=
H^0\bigl(X',(\Omega_{X'}^1)^{\otimes2}\bigr), \qquad F^*a_2(\nabla_t)=\det\bigl(\psi_{\nabla_t}\bigr).
\]
This descent is algebraic in families \cite[Proposition~3.2]{LP01}, and defines 
\[
\begin{aligned}
	h_{\Hod, V}\colon\mathrm{tConn}^0(V) &\longrightarrow \mathcal A_{X'}^0\times \Bbbk,\qquad \nabla_t \longmapsto \bigl(a_2(\nabla_t),t\bigr),
\end{aligned}% = \bigl(\operatorname{Char}(\psi_{\nabla_t})^{1/p},t\bigr).
\]
which is called the \textit{restricted Hodge--Hitchin morphism} to $\mathrm{tConn}^0(V)$ \cite[Section 2.2]{dCZ22}, where $F^*a_2(\nabla_t)=\det\bigl(\psi_{\nabla_t}\bigr)$.

We now use the vector-space structure and coordinates $x_1,\ldots,x_n,t$ introduced above.
\begin{prop}\label{prop:t-Hitchin-homogeneity}
	The morphism $h_{\Hod, V}$ is $\mathbb G_m$-equivariant for scalar multiplication on $\mathrm{tConn}^0(V)$. More precisely, the $\mathbb G_m$-action on $\mathcal A_{X'}^0\times \Bbbk$ is given by
	\[
	\lambda\cdot(a_2,t) = (\lambda^{2p}a_2,\lambda t)
	\]
    In particular, its coordinate functions are homogeneous of degrees $2p$ and $1$, respectively.
\end{prop}

\begin{proof}
	For $\lambda\in\Bbbk$, scalar multiplication sends $\nabla_t$ to the determinant-normalized $(\lambda t)$-connection $\lambda\nabla_t$. For every local vector field $D$, the definition of $p$-curvature gives
	\[
	\begin{aligned}
		\psi(\lambda\nabla_t)(D) &=(\lambda\nabla_t(D))^p - (\lambda t)^{p-1}\lambda\nabla_t(D^{[p]})\\
		&=\lambda^p \bigl(\nabla_t(D)^p - t^{p-1}\nabla_t(D^{[p]})\bigr)\\
		&=\lambda^p\psi_{\nabla_t}(D).
	\end{aligned}
	\]
	Taking determinants, we obtain
	\[
	\det\bigl(\psi(\lambda\nabla_t)\bigr) = \lambda^{2p}\det\bigl(\psi_{\nabla_t}\bigr).
	\]
	Consequently,
	\[
	F^*a_2(\lambda\nabla_t)	= F^*\bigl(\lambda^{2p}a_2(\nabla_t)\bigr).
	\]
	Since pullback under the relative Frobenius is injective and $\Bbbk$-linear on these spaces of sections, it follows that
	\[
	a_2(\lambda\nabla_t) = \lambda^{2p}a_2(\nabla_t).
	\]
	Together with the scaling of $t$, this proves the assertion.
\end{proof}

Combining Proposition \ref{prop:t-Hitchin-homogeneity} with \cite[Theorem 3.5]{Ze18}, we obtain the following result:
\begin{coro}\label{cor:very-stable-nilpotent-connections}
	Let $V$ be a very stable vector bundle on $X$. Then, the restricted Hodge--Hitchin morphism $h_{\Hod, V}$ is finite, flat and surjective. In particular, the set of ordinary connections $\nabla$ on $V$ satisfying $\det\nabla=d$ and having nilpotent $p$-curvature is finite and nonempty.
\end{coro}

\begin{proof}
	A point of $h_{\Hod, V}^{-1}(0,0)$ has $t=0$ and therefore corresponds to a trace-free Higgs field $\phi$. The very stability then forces $\phi=0$, so
	\[
	h_{\Hod, V}^{-1}(0,0)=\{0\}.
	\]	
	  Because of \eqref{eq:dim-tConn}, we see that both source and target are vector spaces of dimension
	\[
	\dim\mathrm{tConn}^0(V)	= 3g-2 = \dim\bigl(\mathcal A_{X'}^0 \times \Bbbk\bigr).
	\]
	By Proposition~\ref{prop:t-Hitchin-homogeneity}, the coordinate functions of the restricted Hodge--Hitchin morphism $h_{\Hod, V}$ are homogeneous of positive degrees. The homogeneous finiteness criterion \cite[Theorem~3.5]{Ze18} therefore implies that $h_{\Hod, V}$ is finite.

	The image of this finite morphism is closed and has the same dimension as its source. Since the target is irreducible of the same dimension, so the restricted Hodge--Hitchin morphism $h_{\Hod, V}$ is surjective. In particular, $h_{\Hod, V}^{-1}(0,1)$ is finite and nonempty. Observe that its points are precisely the ordinary connections on $V$ with determinant connection $d$ and nilpotent $p$-curvature.
	
	Finally, the source is Cohen--Macaulay and the target is regular, since both are affine spaces. Every fibre is zero-dimensional by finiteness, and this equals the difference of their dimensions. The morphism is flat by the miracle flatness theorem \cite[Theorem~23.1]{Ma86}. Therefore, $h_{\Hod, V}$ is flat.
\end{proof}

\subsection{Finiteness of the restricted de Rham--Hitchin morphism}

In \cite[Theorem 4.15]{dCGZ24} de Cataldo, Groechenig, and Zhang establish a Nonabelian Hodge correspondence using very good splittings. A choice of such a splitting identifies the corresponding Higgs and flat bundles, preserving stability and semistability.

Let \(S^{\mathrm{vg}}_{\mathrm{GL}_2}\to\mathcal A_{X'}\) be the moduli scheme of the stack of very good splittings defined in \cite[Definition 3.5, Theorem 4.13]{dCGZ24}. For the inclusion 
\(\mathcal A_{X'}^0\hookrightarrow\mathcal A_{X'} := H^0(X',\Omega_{X'}) \oplus \mathcal A_{X'}^0 \), set
\[
S^{\mathrm{vg}}_{\mathrm{SL}_2}:=
S^{\mathrm{vg}}_{\mathrm{GL}_2}\times_{\mathcal A_{X'}}\mathcal A_{X'}^0.
\]

In the present nilpotent setting, the required very good splitting exists on the double curve \(2X'\subset T^*X'\), which is the spectral curve over the origin. Indeed, the canonical splitting \(F_*\mathcal O_X\) on the reduced zero section extends to a very good splitting on \(2X'\), as shown in \cite[Lemma 3.16]{dCGZ24}. 
This also gives a point of \(S^{\mathrm{vg}}_{\mathrm{SL}_2}\) over the origin. Applying \cite[Theorem 4.15]{dCGZ24} after this base change, with the compatible determinant conditions, yields the following correspondence for the \(\mathrm{SL}_2\)-nilpotent fibres. 
\begin{lemm}\label{lem:scheme-IC}
	The inverse Cartier transform induces an isomorphism of scheme-theoretic nilpotent cones
	\[
	\mathrm{OV}_h:  h^{-1}_{\Dol}(0) \stackrel{\sim}{\lra}   h^{-1}_{\dR}(0)
	\]
	It sends $C'_k$ to $C_k$  for every $k$ and preserves their multiplicities. \qed
\end{lemm}

%\begin{proof}By base change, this lemma is an immediate consequence of \cite[Theorem 4.15]{dCGZ24}.\end{proof}

Now we decompose the Dolbeault nilpotent cone $\mathcal{N}_{\Dol}$ into a union of equidimensional components. According to \eqref{eq:Dol-nilpotent-components} below, we have
\[
h_{\Dol}^{-1}(0)_{red}= \mathcal{N}_{\Dol} = C'_0 \cup \bigcup_{k=2}^g C'_k,
\]
where $C'_0\simeq\mathcal M(X',2)$ is the zero-Higgs component. We would like to call it by \textit{the Cartier component}. For each $k\ge 0$, denote by $C_k:=\mathrm{OV}_h(C'_k)$. 
Then, $C_k$ for $k\ne g$ are precisely the distinct irreducible components of $\mathcal N_{\dR}$, and the equidimensionality of the nilpotent cone \cite{La88} gives
\[
h_{\dR}^{-1}(0)_{red}=\mathcal N_{\dR} = C_0 \cup \bigcup_{k=2}^g C_k, \qquad \dim C_k=\dim C'_k=3g-3.
\]

On $C_0$, the image $\mathrm{OV}_h(E,0)=(F^*E,\nabla^{\rm can})$ is independent of the class $h$ \cite[Theorem~2.8]{OV07} and \cite[Section~2.2]{LSZ15}. Thus, the de Rham forgetful map on $C_0$ identifies with the rational Frobenius-pullback map, called the generalized Verschiebung
\[
Ver: \mathcal M(X',2)\dashrightarrow\mathcal M(X,2), \qquad [E]\longmapsto[F^*E],
\]
defined where $F^*E$ is semistable. 

For $k\ge 2$, the general point of $C'_k$ has nonzero Higgs field, and the general point of $C_k$ has nonzero nilpotent $p$-curvature. Accordingly, we write
\[
\mathcal N_{\dR}^{\ne0}:= \bigcup_{k=2}^g C_k
\]
for the union of these components. Notice that this union may meet $C_0$. We refer to the rational map
$$ \pi_{\dR}: \mathcal{N}_{\dR}^{\not= 0} \dashrightarrow \mathcal{M}(X,2), \qquad (V,\nabla) \mapsto V$$
as the \textit{rational de Rham forgetful map}.

Let $\mathcal M_X^{\vs}\subset\mathcal M(X,2)$ be the open subset of very stable bundles, and fix $V\in\mathcal M_X^{\vs}$. Choose an ordinary connection $\nabla_0$ with $\det\nabla_0=d$, as in the preceding subsections. Denote by $W_V:=H^0(X,\End_0(V)\otimes\Omega_X^1)$ and
\[
\Conn^0(V):=\left\{ \nabla\colon V\longrightarrow V\otimes\Omega_X^1 \ \middle|\ \nabla\text{ is a connection and }\det\nabla=d \right\}.
\]
Then
\[
\Conn^0(V)=\nabla_0+W_V, \qquad \dim W_V=\dim\Conn^0(V)=n:=3g-3.
\]

Recall that
\[
\mathcal A_{X'}^0 = H^0(X',\Omega_{X'}^{\otimes2}).
\]
The de Rham--Hitchin morphism of \cite[Proposition~3.2]{LP01}, restricted to connections on $V$ with trivial determinant connection, is
\[
h_{\dR,V}\colon \Conn^0(V)\longrightarrow\mathcal A_{X'}^0, \qquad \nabla\longmapsto a_2(\nabla),
\]
where $a_2(\nabla)$ is uniquely determined by
\[
F^*a_2(\nabla)=\det(\psi_\nabla).
\]
Equivalently, $h_{\dR,V}$ is obtained from $h_{\Hod, V}$ by base change along
\[
\mathcal A_{X'}^0\longrightarrow \mathcal A_{X'}^0\times\Bbbk, \qquad a\longmapsto(a,1).
\]

\begin{prop}\label{prop:restricted}
	For every $V\in\mathcal M_X^{\vs}$, the restricted de Rham--Hitchin morphism
	\[
	h_{\dR,V}\colon \Conn^0(V)\longrightarrow\mathcal A_{X'}^0
	\]
	is finite, flat and surjective, of degree $\deg(h_{\dR,V}) = (2p)^{3g-3}$. In particular, every geometric fibre has scheme-theoretic length $(2p)^{3g-3}$.
\end{prop}

\begin{proof}
	By Corollary~\ref{cor:very-stable-nilpotent-connections}, the morphism $h_{\Hod, V}$ is finite, flat and surjective. These properties are preserved by base change, so they hold for $h_{\dR,V}$. To compute its degree, choose bases $\phi_1,\ldots,\phi_n$ of $W_V$ and $\eta'_1,\ldots,\eta'_n$ of $\mathcal A_{X'}^0$. Write
	\[
	\phi_x=\sum_{i=1}^n x_i\phi_i, \qquad a_2(t\nabla_0+\phi_x) = \sum_{j=1}^n Q_j(x,t)\eta'_j.
	\]
	By Proposition~\ref{prop:t-Hitchin-homogeneity}, each $Q_j$ is homogeneous of total degree $2p$. Set
	\[
	q_j(x):=Q_j(x,0), \qquad R:=\Bbbk[x_1,\ldots,x_n].
	\]
	The fibre of $h_{\Hod, V}$ over $(0,0)$ is therefore
	\[
	h_{\Hod, V}^{-1}(0,0) = \operatorname{Spec}R/(q_1,\ldots,q_n).
	\]
	
	At $t=0$, the $0$-connection is the Higgs field $\phi_x$, and
	\[
	F^*a_2(\phi_x) = (\det\phi_x)^{\otimes p}.
	\]
	Notice that the condition $a_2(\phi_x)=0$ is equivalent to nilpotence of $\phi_x$. Then, by using the Nullstellensatz theorem, very stability of $V$ implies that
	\[
	\sqrt{(q_1,\ldots,q_n)} = (x_1,\ldots,x_n).
	\]
	Consequently, the $q_j$ form a homogeneous system of parameters in the Cohen--Macaulay ring $R$, and hence a regular sequence. Observe that each $q_i$ has degree $2p$, so
	\[
	\operatorname{Hilb}_{R/(q_1,\ldots,q_n)}(z)	= \frac{(1-z^{2p})^n}{(1-z)^{3g-3}}
	=
	(1+z+\cdots+z^{2p-1})^{3g-3}.
	\]
	Evaluating at $z=1$ yields
	\[
	\operatorname{length} \bigl(h_{\Hod, V}^{-1}(0,0)\bigr) = (2p)^{3g-3}.
	\]
	
	Finally, $h_{\Hod, V}$ is finite and flat over the connected base $\mathcal A_{X'}^0\times\Bbbk$. Its direct image of the structure sheaf therefore has constant rank, equal to this fibre length. The base change to $t=1$ has the same rank, proving
	\[
	\deg(h_{\dR,V})=(2p)^{3g-3}.
	\]
\end{proof}

\section{Geometry of the rank-$2$ de Rham nilpotent cone}
\subsection{The Dolbeault and de Rham nilpotent cones}

Fix $h\in\LL(X)$ and assume that $g\ge2$ and $p \not= 2$. We recall the degree-indexed description of the rank-$2$ Dolbeault nilpotent cone  $\mathcal{N}_{Dol}$ with trivial determinant, following \cite[Section~2, especially Proposition~2.2]{PP21}. The zero-Higgs component is
\[
C'_0 := \left\{[(E,0)]\mid[E]\in\mathcal M(X',2)\right\} \simeq\mathcal M(X',2).
\]

Let $(E,\phi)$ be a semistable Higgs bundle with nonzero nilpotent Higgs field. Its extension data are
\[
0\longrightarrow L\xrightarrow{i_E}E \xrightarrow{\pi_E}L^{-1}\longrightarrow0, \qquad e\in H^1(X',L^2), \qquad 0\ne u\in H^0(X',L^2\Omega_{X'}),
\]
where $L=\ker\phi$ and $\phi=(i_E\otimes\id_{\Omega_{X'}})\circ u\circ \pi_E$. Following the indexing of Pal--Pauly, write $\deg L=1-k$. As in \cite[Section~2]{PP21}, semistability of $(E,\phi)$ and $u\ne0$ give $1\le k\le g$.

For $2\le k<g$, let $C'_k$ be the closure of the locus of such Higgs bundles with $\deg(\ker\phi)=1-k$. The parameterization in \cite[Proposition~2.2]{PP21}, together with the discussion immediately following its proof, identifies $C'_k$ as an irreducible component of dimension
\[
\dim C'_k = (2g-2k)+(g+2k-3) = 3g-3.
\]

The case $k=g$ requires a separate description. In this case, $\deg L=1-g$, so $L^2\Omega_{X'}$ has degree zero. A nonzero section of a degree-zero line bundle has no zeros and therefore trivializes that bundle. Consequently, the condition $0\ne u\in H^0(X',L^2\Omega_{X'})$ implies that $L^2\Omega_{X'}\simeq\OO_{X'}$. Hence, the line bundle $\theta:=L^{-1}$ satisfies $\theta^{\otimes2}\simeq\Omega_{X'}$. This shows that $\theta$ is a theta characteristic, determined by the Higgs pair through $\theta=(\ker\phi)^{-1}$. This is precisely the description of the possible line bundles $L$ in \cite[Section~5.3]{PP21}.

Denote the set of theta characteristics by
\[
\Theta(X'):= \left\{ \theta \in\Pic^{g-1}(X') \;\middle|\; \theta^{\otimes2}\simeq\Omega_{X'}\right\}.
\]
It is a torsor under $\Pic^0(X')[2]$ and has $2^{2g}$ elements, since $p>2$. For each $\theta\in\Theta(X')$, define \(C'_{g,\theta}\) to be the closure of the locus of Higgs bundles \((E,\phi)\) constructed from an extension
\[
0\longrightarrow\theta^{-1} \xrightarrow{i_E}E \xrightarrow{\pi_E}\theta \longrightarrow0,
\]
with extension class \(e\in H^1(X',\theta^{-2})\), together with a nonzero morphism
\[
u\colon\theta\longrightarrow\theta^{-1}\otimes\Omega_{X'}.
\]
Proposition 2.2 of \cite{PP21} gives a generically injective rational parameterization of $C'_{g,\theta}$ by
\[
\PP\!\left( H^1(X', T_{X'})\oplus\Bbbk \right) \simeq\PP^{3g-3}.
\]
In particular, $C'_{g,\theta}$ is irreducible of dimension $3g-3$. Thus we have the decomposition 
\[
C'_g = \bigcup_{\theta\in\Theta(X')}C'_{g,\theta}.
\]
The case $k=1$ contributes no additional irreducible component. Indeed, $\deg L=0$, and the Higgs Jordan--Hölder filtration
\[
0\subset(L,0)\subset(E,\phi)
\]
gives
\[
(E,\phi)\sim_S(L\oplus L^{-1},0).
\]
Its image in the coarse moduli space is therefore
\begin{equation}\label{eq:Kummer-component}
    C'_1=\mathrm{Kum}_{X'}\subset C'_0.
\end{equation}
Here, the Kummer variety is defined by
\[
\mathrm{Kum}_{X'} := \Pic^0(X')/\iota =  \left\{ [M\oplus M^{-1}] \;\middle|\; M\in\Pic^0(X') \right\} \subset\mathcal M(X',2),
\]
where $\iota$ denotes the involution $\iota(M) = M^{-1}$. We consequently obtain the decomposition into irreducible components
\begin{equation}\label{eq:Dol-nilpotent-components}
	\mathcal N_{\Dol} = h_{Dol}^{-1}(0)_{red} =  C'_0 \cup \bigcup_{k=2}^{g-1}C'_k \cup \bigcup_{\theta\in\Theta(X')}C'_{g,\theta}.
\end{equation}

Fix $h\in\LL(X)$ and set
\[
C_k:=\mathrm{OV}_h(C'_k), \qquad C_{g,\theta}:=\mathrm{OV}_h(C'_{g,\theta}).
\]
By Lemma~\ref{lem:scheme-IC}, the inverse Cartier transform identifies the scheme-theoretic nilpotent fibres and preserves the multiplicities of corresponding irreducible components. It therefore transports \eqref{eq:Dol-nilpotent-components} to
\begin{equation} \label{eq:dR-nilpotent-components}
\mathcal N_{\dR} = h^{-1}_{\dR}(0)_{red} = C_0 \cup \bigcup_{k=2}^{g-1}C_k \cup \bigcup_{\theta\in\Theta(X')}C_{g,\theta}, \qquad C_g = \bigcup_{\theta\in\Theta(X')}C_{g,\theta}.
\end{equation}
In particular, all these de Rham components also have dimension $3g-3$.

On the zero-Higgs component,
\[
\mathrm{OV}_h(E,0) = (F^*E,\nabla^{\rm can}).
\]
For nonzero Higgs extension data, Proposition~\ref{prop:inverse-class} gives
\[
0\longrightarrow F^*L\longrightarrow V \longrightarrow F^*L^{-1}\longrightarrow0, \qquad [V]=F^*e+h\smile F^*u.
\]
The induced connections on the subbundle and quotient are canonical, and
\[
(\det V,\det\nabla)\simeq(\OO_X,d).
\]
For $C_{g,\theta}$, this specializes to
\[
0\longrightarrow F^*\theta^{-1} \longrightarrow V \longrightarrow F^*\theta \longrightarrow0.
\]

%\begin{rema}[Proof of Proposition \ref{dominance}(1)] \label{rem:degree-zero-Kummer}
%	For $k=1$, both line bundles in the de Rham extension have degree zero. Thus $V$ is strictly semistable, and $0\subset F^*L\subset V$ is a Jordan--%Hölder filtration of its underlying vector bundle, with
%	\[
%	\operatorname{gr}(V) \simeq F^*L\oplus F^*L^{-1}.
%	\]
%	Hence
%	\[
%	V\sim_S F^*L\oplus F^*L^{-1},
%	\]
%	where $\sim_S$ denotes S-equivalence, namely isomorphism of the Jordan--Hölder graded bundles. The relation expresses equality of their points in  % %$\mathcal M(X,2)$; the extension need not split. Consequently, the forgetful image of $C_1$ lies in $\mathrm{Kum}_X$ and is not dominant in
%$\mathcal M(X,2)$. This show that the restricted map on $C_1$ is not dominant.
%\end{rema}

%With the indexing $\deg L=1-k$, the Hausel--Hitchin multiplicity formula for the corresponding complex type-$(1,1)$ components is
%\[
%m_k=2^{2g-2k}, \qquad 2\le k\le g;
%\]
%see \cite[Corollary~5.15]{HH2}. Lemma~\ref{lem:scheme-IC} identifies the multiplicities of corresponding Dolbeault and de Rham components in our setting.

\subsection{The dominance of the rational de Rham forgetful map}

In this subsection we show Proposition \ref{dominance}. Let $V$ be a very stable bundle in $\mathcal{M}(X,2)$. We know by Proposition \ref{prop:restricted} that the map $h_{dR,V}$ is finite, flat of degree $(2p)^{3g-3}$. Thus the fibre $h_{dR,V}^{-1}(0)$ is a $0$-dimensional scheme of length $(2p)^{3g-3}$. On the other hand, if we denote by $m_k$ the multiplicity of the reduced component $C_k$ appearing in the decomposition \eqref{eq:dR-nilpotent-components} of $\mathcal{N}_{\dR}$ in the scheme-theoretic nilpotent cone $h^{-1}_{\dR}(0)$, then we obtain a decomposition of $0$-dimensional schemes
$$h_{dR,V}^{-1}(0) = m_0 (C_0 \cap\mathrm{Conn^0(V)}) \cup  \bigcup_{k=2}^g m_k (C_k \cap\mathrm{Conn^0(V)}). $$
Considering the restricted Dolbeault--Hitchin morphism for a very stable bundle and using Lemma \ref{lem:scheme-IC} it is easily seen that $m_0 = 2^{3g-3}$.
Also, the length of the scheme $C_0 \cap \mathrm{Conn^0(V)}$ coincides with the degree of the generalized Verschiebung $\deg (Ver)$ --- see e.g. 
\cite{Os06}. Therefore, taking the length of the schemes in the previous decomposition, we obtain the equality
$$ (2p)^{3g-3} = 2^{3g-3} \deg(Ver) + res,$$
where $res$ denotes the residual length. It will therefore be enough to show that $res > 0$, or equivalently
$\deg(Ver) < p^{3g-3}$.

\bigskip
We know that the generalized Verschiebung 
$$ Ver : \mathcal{M}(X',2) \dashrightarrow \mathcal{M}(X,2) , \qquad [E] \mapsto [F^*E]$$
is a dominant \cite[Theorem A.6]{Os06}, rational map. The associated degree of the function field extension equals $p^{3g-3}$, which follows
from the fact that $Ver^* (\mathcal L) = {\mathcal L}'^p$, where $\mathcal{L}$ and $\mathcal{L}'$ denote the determinant line bundles over the
moduli spaces \cite[Proposition 7.1]{LP04}. Moreover, the indeterminacy locus of $Ver$ is non-empty for any curve (see e.g. 
\cite[Proposition 7.4]{LP04}), which shows the strict inequality $\deg(Ver) < p^{3g-3}$.

\begin{rema}
An explicit formula for the degree of the generalized Verschiebung is known thanks to \cite{KW25} and \cite{Zh26}, but the formula
only holds for general curves.
\end{rema}

\begin{rema}
	For $g=2$, we have $\mathcal{N}_{\dR} = C_0 \cup C_2$ with forgetful maps
	\[
	C_0\dashrightarrow \mathcal M(X,2),
	\qquad
	C_2\dashrightarrow \mathcal M(X,2),
	\]
	 both being dominant.
\end{rema}

\section{Open questions}

Let $L$ (resp. $L'$) be a line bundle over $X$ (resp. $X'$). 
We denote by $\mathcal{M}(X,2,L)$ the coarse moduli space of semi-stable rank-$2$ vector bundles with
fixed determinant $L$ and similarly for the moduli space $\mathcal{M}(X',2,L')$.  For simplicity we denote 
$\mathcal{M}(X,2) = \mathcal{M}(X,2,\mathcal{O}_X)$. 
We consider the rational map, defined on the non-empty (see section 4.3) open subset where the rank-$2$ vector bundle $\mathcal{H}_{\hat{h}}$ is semi-stable \footnote{Note the equality of $S$-equivalence classes $[\HH_{\hat{h}}]=[\HH_{\lambda\hat{h}}]$ for any 
$\lambda \in \mathbb{G}_m$.}
\begin{equation}\label{eq:Phi}
 \overline \Phi\colon \PP(\widehat{\LL(X)}) = \PP^{3g-3}
 \dashrightarrow \mathcal{M}(X,2,F^*\Omega_{X'}),\qquad
 [\hat{h}]\longmapsto[\HH_{\hat{h}}].
\end{equation}

Note the set-theoretical decomposition $\PP (\widehat{\LL(X)}) = \LL(X) \sqcup 
\PP H^1(X', T_{X'})$. The restriction of $\overline \Phi$ to
 $\PP H^1(X', T_{X'})$ factorizes as follows
 $$ \PP^{3g-2} = \PP H^1(X', T_{X'}) \dashrightarrow \Theta \subset \mathcal{M}(X',2,\Omega_{X'}) \stackrel{Ver}\dashrightarrow 
  \mathcal{M}(X,2,F^* \Omega_{X'}),  $$
 where $\Theta = \{ [E]\ | \ \dim H^0(X',E) \not=0 \}$ is the generalized Theta divisor and $Ver$ denotes the generalized
 Verschiebung $Ver(E) = F^*E$. It can be shown that the closure $\overline{Ver(\Theta)}$ is a divisor. We therefore conjecture

\begin{conj}\label{conj:dominance}
The map \eqref{eq:Phi} is dominant.
\end{conj}

Fix $h\in\LL(X)$ and a theta characteristic $\theta$ on $X'$.
The extension class $e\in H^1(X',T')$ determines a nilpotent Higgs bundle
\[
 0\longrightarrow \theta^{-1} \longrightarrow E_e
 \longrightarrow \theta\longrightarrow 0.
\]
with nilpotent Higgs field given by $t\cdot1  \colon \theta \to \theta^{-1}\Omega_{X'} = \theta$, for $t\in\Bbbk$.
If we denote by  $V_{e,t}$  the underlying bundle of
$\mathrm{OV}_h(E_e,t\cdot1)$, then by Proposition \ref{prop:inverse-class}
\begin{equation}\label{eq:theta-family}
 [V_{e,t}]=F^*e+th \in \widehat{\mathcal{L}(X)} \subset H^1(X,F^* T_{X'}). 
 \end{equation}
Thus, after tensoring with $F^* \theta$, which changes the determinant from $\mathcal{O}_X$ to $F^* \Omega_{X'}$, we obtain that 
\begin{equation} \label{eq:theta-family2}
V_{e,t}\otimes F^*\theta\simeq\HH_{F^*e+th}.
\end{equation}

\bigskip

We have a decomposition $C_{g} = \bigcup_\theta C_{g,\theta}$ where $\theta$ varies over all theta-characteristics of $X'$. 
The previous observation shows that $C_{g,\theta}$ parameterizes flat bundles of the form \eqref{eq:theta-family}.
By restriction of $\pi_{\dR}$ to $C_{g,\theta}$ we obtain a rational forgetful map
\[
 \pi_{g,\theta}\colon C_{g,\theta}
 \dashrightarrow \mathcal{M}(X,2).
\]

\begin{prop}\label{prop:dominance-equivalence}
Conjecture \ref{conj:dominance} is equivalent to dominance of
$\pi_{g,\theta}$ for any theta characteristic $\theta$.
\end{prop}
\begin{proof}
Every $\hat{h} \in \widehat{\LL(X)}$ can be uniquely written as
$\hat{h} = F^*e+th$, by \eqref{eq:cohomology-sequence}.
We note that tensor product with $F^*\theta$ induces an isomorphism
$M(X,2) \xrightarrow{\sim} M(X,2,F^*\Omega_{X'})$. It follows then from
equality \eqref{eq:theta-family2} that the images of the two rational maps 
$\overline{\Phi}$ and $\pi_{g,\theta}$ coincide under the above isomorphism.
\end{proof}

%The restriction of $\overline\Phi$ to the boundary in
%\eqref{eq:projective-compactification} is Frobenius pullback on
%the extension family
%$0\to\OO_{X'}\to E_0\to\Omega_{X'}\to0$, wherever $F^*E_0$ is
%semistable. Thus its image is a Verschiebung image of that
%extension locus. If this boundary image is a divisor in the
%target, finding a point of $\im\overline\Phi$ outside it proves
%dominance, since the total image closure is irreducible.
%The dimension assumption on the boundary image is needed for
%this argument.

%\begin{ques}\label{question:very-stable-H}
%Is $\HH_h$ very stable for general $h\in\LL(X)$?
%Can the Prym construction above produce a Frobenius pullback
%boundary point that is very stable as a bundle on $X$?
%\end{ques}
%This asks for more than stability of the flat pair and more than
%very stability before Frobenius pullback. Together with
%Conjecture~\ref{conj:dominance}, it isolates the underlying-bundle
%questions that remain after the extension-class description of
%the correspondence.

A natural question related to Proposition \ref{dominance} is
\begin{ques}
For which integers $k$ is the restriction $\pi_k$ of the rational forgetful map $\pi_{\dR}$ to the component $C_k$ dominant?
\end{ques}

More generally, if we denote by $N_k$ the degree of the restriction $\pi_k$  (and putting $N_k = 0$ if $\pi_k$ is not dominant), then we might expect the following expression for the length of the $0$-dimensional scheme $h^{-1}_{dR,V}(0)$ for a general bundle $V$

\begin{conj}
We have the following relation $(2p)^{3g-3} = 2^{3g-3} N_0 + \sum_{k=2}^g 2^{2g-2k} N_k$.
\end{conj}

Note that $N_0$ is the degree of the generalized Verschiebung computed recently by \cite{KW25} and \cite{Zh26}.
The coefficients in front of the $N_k$ are the multiplicities of the reduced components $C_k$ for the scheme-structure on 
the de Rham nilpotent cone $h^{-1}_{\dR}(0)$. They were given in \cite[Corollary 5.15]{HH22} for the Dolbeault 
nilpotent cone $h^{-1}_{\Dol}(0)$ over the complex numbers.

\end{document}